\documentclass[reqno,12pt,a4paper]{amsart}

\usepackage{mathrsfs}
\usepackage{amssymb}
\usepackage{amsfonts}
\usepackage{latexsym}
\usepackage{amsthm}
\usepackage{graphicx}

\def\lb{\label}

\newcommand{\er}[1]{\textrm{(\ref{#1})}}

\begin{document}


\renewcommand{\theequation}{\arabic{section}.\arabic{equation}}
\theoremstyle{plain}
\newtheorem{theorem}{\bf Theorem}[section]
\newtheorem{lemma}[theorem]{\bf Lemma}
\newtheorem{corollary}[theorem]{\bf Corollary}
\newtheorem{proposition}[theorem]{\bf Proposition}
\newtheorem{definition}[theorem]{\bf Definition}
\newtheorem{example}[theorem]{\bf Example}

\theoremstyle{remark}
\newtheorem*{remarks}{\bf Remarks}
\newtheorem*{remark}{\bf Remark}

\def\a{\alpha}  \def\cA{{\mathcal A}}     \def\bA{{\bf A}}  \def\mA{{\mathscr A}}
\def\b{\beta}   \def\cB{{\mathcal B}}     \def\bB{{\bf B}}  \def\mB{{\mathscr B}}
\def\g{\gamma}  \def\cC{{\mathcal C}}     \def\bC{{\bf C}}  \def\mC{{\mathscr C}}
\def\G{\Gamma}  \def\cD{{\mathcal D}}     \def\bD{{\bf D}}  \def\mD{{\mathscr D}}
\def\d{\delta}  \def\cE{{\mathcal E}}     \def\bE{{\bf E}}  \def\mE{{\mathscr E}}
\def\D{\Delta}  \def\cF{{\mathcal F}}     \def\bF{{\bf F}}  \def\mF{{\mathscr F}}
\def\c{\chi}    \def\cG{{\mathcal G}}     \def\bG{{\bf G}}  \def\mG{{\mathscr G}}
\def\z{\zeta}   \def\cH{{\mathcal H}}     \def\bH{{\bf H}}  \def\mH{{\mathscr H}}
\def\e{\eta}    \def\cI{{\mathcal I}}     \def\bI{{\bf I}}  \def\mI{{\mathscr I}}
\def\p{\psi}    \def\cJ{{\mathcal J}}     \def\bJ{{\bf J}}  \def\mJ{{\mathscr J}}
\def\vT{\Theta} \def\cK{{\mathcal K}}     \def\bK{{\bf K}}  \def\mK{{\mathscr K}}
\def\k{\kappa}  \def\cL{{\mathcal L}}     \def\bL{{\bf L}}  \def\mL{{\mathscr L}}
\def\l{\lambda} \def\cM{{\mathcal M}}     \def\bM{{\bf M}}  \def\mM{{\mathscr M}}
\def\L{\Lambda} \def\cN{{\mathcal N}}     \def\bN{{\bf N}}  \def\mN{{\mathscr N}}
\def\m{\mu}     \def\cO{{\mathcal O}}     \def\bO{{\bf O}}  \def\mO{{\mathscr O}}
\def\n{\nu}     \def\cP{{\mathcal P}}     \def\bP{{\bf P}}  \def\mP{{\mathscr P}}
\def\r{\rho}    \def\cQ{{\mathcal Q}}     \def\bQ{{\bf Q}}  \def\mQ{{\mathscr Q}}
\def\s{\sigma}  \def\cR{{\mathcal R}}     \def\bR{{\bf R}}  \def\mR{{\mathscr R}}
\def\S{\Sigma}  \def\cS{{\mathcal S}}     \def\bS{{\bf S}}  \def\mS{{\mathscr S}}
\def\t{\tau}    \def\cT{{\mathcal T}}     \def\bT{{\bf T}}  \def\mT{{\mathscr T}}
\def\f{\phi}    \def\cU{{\mathcal U}}     \def\bU{{\bf U}}  \def\mU{{\mathscr U}}
\def\F{\Phi}    \def\cV{{\mathcal V}}     \def\bV{{\bf V}}  \def\mV{{\mathscr V}}
\def\P{\Psi}    \def\cW{{\mathcal W}}     \def\bW{{\bf W}}  \def\mW{{\mathscr W}}
\def\o{\omega}  \def\cX{{\mathcal X}}     \def\bX{{\bf X}}  \def\mX{{\mathscr X}}
\def\x{\xi}     \def\cY{{\mathcal Y}}     \def\bY{{\bf Y}}  \def\mY{{\mathscr Y}}
\def\X{\Xi}     \def\cZ{{\mathcal Z}}     \def\bZ{{\bf Z}}  \def\mZ{{\mathscr Z}}

\def\vr{\varrho}
\def\be{{\bf e}} \def\bc{{\bf c}}
\def\bx{{\bf x}} \def\by{{\bf y}}
\def\bv{{\bf v}} \def\bu{{\bf u}}
\def\Om{\Omega} \def\bp{{\bf p}}
\def\bbD{\pmb \Delta}
\def\mm{\mathrm m}
\def\mn{\mathrm n}

\newcommand{\mc}{\mathscr {c}}

\newcommand{\gA}{\mathfrak{A}}          \newcommand{\ga}{\mathfrak{a}}
\newcommand{\gB}{\mathfrak{B}}          \newcommand{\gb}{\mathfrak{b}}
\newcommand{\gC}{\mathfrak{C}}          \newcommand{\gc}{\mathfrak{c}}
\newcommand{\gD}{\mathfrak{D}}          \newcommand{\gd}{\mathfrak{d}}
\newcommand{\gE}{\mathfrak{E}}
\newcommand{\gF}{\mathfrak{F}}           \newcommand{\gf}{\mathfrak{f}}
\newcommand{\gG}{\mathfrak{G}}           
\newcommand{\gH}{\mathfrak{H}}           \newcommand{\gh}{\mathfrak{h}}
\newcommand{\gI}{\mathfrak{I}}           \newcommand{\gi}{\mathfrak{i}}
\newcommand{\gJ}{\mathfrak{J}}           \newcommand{\gj}{\mathfrak{j}}
\newcommand{\gK}{\mathfrak{K}}            \newcommand{\gk}{\mathfrak{k}}
\newcommand{\gL}{\mathfrak{L}}            \newcommand{\gl}{\mathfrak{l}}
\newcommand{\gM}{\mathfrak{M}}            \newcommand{\gm}{\mathfrak{m}}
\newcommand{\gN}{\mathfrak{N}}            \newcommand{\gn}{\mathfrak{n}}
\newcommand{\gO}{\mathfrak{O}}
\newcommand{\gP}{\mathfrak{P}}             \newcommand{\gp}{\mathfrak{p}}
\newcommand{\gQ}{\mathfrak{Q}}             \newcommand{\gq}{\mathfrak{q}}
\newcommand{\gR}{\mathfrak{R}}             \newcommand{\gr}{\mathfrak{r}}
\newcommand{\gS}{\mathfrak{S}}              \newcommand{\gs}{\mathfrak{s}}
\newcommand{\gT}{\mathfrak{T}}             \newcommand{\gt}{\mathfrak{t}}
\newcommand{\gU}{\mathfrak{U}}             \newcommand{\gu}{\mathfrak{u}}
\newcommand{\gV}{\mathfrak{V}}             \newcommand{\gv}{\mathfrak{v}}
\newcommand{\gW}{\mathfrak{W}}             \newcommand{\gw}{\mathfrak{w}}
\newcommand{\gX}{\mathfrak{X}}               \newcommand{\gx}{\mathfrak{x}}
\newcommand{\gY}{\mathfrak{Y}}              \newcommand{\gy}{\mathfrak{y}}
\newcommand{\gZ}{\mathfrak{Z}}             \newcommand{\gz}{\mathfrak{z}}

\def\ve{\varepsilon}   \def\vt{\vartheta}    \def\vp{\varphi}    \def\vk{\varkappa}

\def\A{{\mathbb A}} \def\B{{\mathbb B}} \def\C{{\mathbb C}}
\def\dD{{\mathbb D}} \def\E{{\mathbb E}} \def\dF{{\mathbb F}} \def\dG{{\mathbb G}} \def\H{{\mathbb H}}\def\I{{\mathbb I}} \def\J{{\mathbb J}} \def\K{{\mathbb K}} \def\dL{{\mathbb L}}\def\M{{\mathbb M}} \def\N{{\mathbb N}} \def\O{{\mathbb O}} \def\dP{{\mathbb P}} \def\R{{\mathbb R}}\def\S{{\mathbb S}} \def\T{{\mathbb T}} \def\U{{\mathbb U}} \def\V{{\mathbb V}}\def\W{{\mathbb W}} \def\X{{\mathbb X}} \def\Y{{\mathbb Y}} \def\Z{{\mathbb Z}}


\def\la{\leftarrow}              \def\ra{\rightarrow}            \def\Ra{\Rightarrow}
\def\ua{\uparrow}                \def\da{\downarrow}
\def\lra{\leftrightarrow}        \def\Lra{\Leftrightarrow}


\def\lt{\biggl}                  \def\rt{\biggr}
\def\ol{\overline}               \def\wt{\widetilde}
\def\ul{\underline}
\def\no{\noindent}


\let\ge\geqslant                 \let\le\leqslant
\def\lan{\langle}                \def\ran{\rangle}
\def\/{\over}                    \def\iy{\infty}
\def\sm{\setminus}               \def\es{\emptyset}
\def\ss{\subset}                 \def\ts{\times}
\def\pa{\partial}                \def\os{\oplus}
\def\om{\ominus}                 \def\ev{\equiv}
\def\iint{\int\!\!\!\int}        \def\iintt{\mathop{\int\!\!\int\!\!\dots\!\!\int}\limits}
\def\el2{\ell^{\,2}}             \def\1{1\!\!1}
\def\sh{\sharp}
\def\wh{\widehat}
\def\bs{\backslash}
\def\intl{\int\limits}

\def\na{\mathop{\mathrm{\nabla}}\nolimits}
\def\sh{\mathop{\mathrm{sh}}\nolimits}
\def\ch{\mathop{\mathrm{ch}}\nolimits}
\def\where{\mathop{\mathrm{where}}\nolimits}
\def\all{\mathop{\mathrm{all}}\nolimits}
\def\as{\mathop{\mathrm{as}}\nolimits}
\def\Area{\mathop{\mathrm{Area}}\nolimits}
\def\arg{\mathop{\mathrm{arg}}\nolimits}
\def\const{\mathop{\mathrm{const}}\nolimits}
\def\det{\mathop{\mathrm{det}}\nolimits}
\def\diag{\mathop{\mathrm{diag}}\nolimits}
\def\diam{\mathop{\mathrm{diam}}\nolimits}
\def\dim{\mathop{\mathrm{dim}}\nolimits}
\def\dist{\mathop{\mathrm{dist}}\nolimits}
\def\Im{\mathop{\mathrm{Im}}\nolimits}
\def\Iso{\mathop{\mathrm{Iso}}\nolimits}
\def\Ker{\mathop{\mathrm{Ker}}\nolimits}
\def\Lip{\mathop{\mathrm{Lip}}\nolimits}
\def\rank{\mathop{\mathrm{rank}}\limits}
\def\Ran{\mathop{\mathrm{Ran}}\nolimits}
\def\Re{\mathop{\mathrm{Re}}\nolimits}
\def\Res{\mathop{\mathrm{Res}}\nolimits}
\def\res{\mathop{\mathrm{res}}\limits}
\def\sign{\mathop{\mathrm{sign}}\nolimits}
\def\span{\mathop{\mathrm{span}}\nolimits}
\def\supp{\mathop{\mathrm{supp}}\nolimits}
\def\Tr{\mathop{\mathrm{Tr}}\nolimits}
\def\BBox{\hspace{1mm}\vrule height6pt width5.5pt depth0pt \hspace{6pt}}


\newcommand\nh[2]{\widehat{#1}\vphantom{#1}^{(#2)}}
\def\dia{\diamond}

\def\Oplus{\bigoplus\nolimits}



\def\qqq{\qquad}
\def\qq{\quad}
\let\ge\geqslant
\let\le\leqslant
\let\geq\geqslant
\let\leq\leqslant
\newcommand{\ca}{\begin{cases}}
\newcommand{\ac}{\end{cases}}
\newcommand{\ma}{\begin{pmatrix}}
\newcommand{\am}{\end{pmatrix}}
\renewcommand{\[}{\begin{equation}}
\renewcommand{\]}{\end{equation}}
\def\eq{\begin{equation}}
\def\qe{\end{equation}}
\def\[{\begin{equation}}
\def\bu{\bullet}

\title[Trace formulas for Schr\"odinger operators on periodic graphs]
{Trace formulas for Schr\"odinger operators on periodic graphs}

\date{\today}
\author[Evgeny Korotyaev]{Evgeny Korotyaev}
\address{Saint-Petersburg State University, Universitetskaya nab. 7/9, St. Petersburg, 199034, Russia,
and HSE University, 3A Kantemirovskaya ulitsa, St. Petersburg, 194100,
Russia, \ korotyaev@gmail.com, \
e.korotyaev@spbu.ru,}
\author[Natalia Saburova]{Natalia Saburova}
\address{Northern (Arctic) Federal University, Severnaya Dvina Emb. 17, Arkhangelsk, 163002, Russia,
 \ n.saburova@gmail.com, \ n.saburova@narfu.ru}


\begin{abstract}
We consider Schr\"odinger operators with periodic potentials on
periodic discrete graphs. Their spectrum consists of a finite number
of bands. We determine trace formulas for the Schr\"odinger
operators. The proof is based on the decomposition of the
Schr\"odinger operators into a direct integral and a specific representation of fiber operators. The traces of the fiber operators are expressed as finite Fourier series of the quasimomentum. The coefficients of the Fourier series are given in terms of the potentials and cycles in the quotient graph from some specific cycle sets. We also present the trace formulas for the heat kernel and the resolvent of the Schr\"odinger operators and the determinant formulas.
\end{abstract}

\maketitle

\section {\lb{Sec1}Introduction}
\setcounter{equation}{0}

There are a lot of applications of different periodic media, e.g.
nanomedia, in physics, chemistry and engineering, see, e.g.,
\cite{NG04}. In order to study properties of such media
one uses their approximations by periodic graphs. In these cases we need to investigate various properties of Laplace and Schr\"odinger operators on
graphs. It is known that the spectrum of Schr\"odinger operators with periodic
potentials on periodic discrete graphs consists of an absolutely
continuous part (a~union of a finite number of non-degenerate bands)
and a finite number of flat bands, i.e., eigenvalues  of infinite
multiplicity. One of the interesting problems is trace formulas,
which show the relationship between the spectra (or spectral shift
functions) of these operators and geometric parameters of graphs.

The trace formulas for multidimensional Schr\"odinger operators on
the lattice $\Z^d$ with real decaying potentials  were obtained in
\cite{IK12}. Trace formulas for Schr\"odinger operators with complex decaying
potentials $V$ on the lattice $\Z^d$ were determined in \cite{KL18},
\cite{K17} and for the case $\Im V\le 0$ in \cite{MN15}. Trace
formulas for time periodic complex Hamiltonians on lattice were
determined in \cite{K21}.

Trace formulas for discrete Laplacians on finite (mostly regular)
graphs have been developed by various authors, see \cite{Ah87,Br91,Mn07,OGS09}. The traces for the $n$-th power of the Laplacian were expressed as the numbers of cycles of length $n$ in the graph.

Trace formulas for discrete Laplacians on infinite graphs are much less understood. It is not clear how to relate the spectrum of the Laplacian on an arbitrary infinite graph with the number of cycles in the graph, since in this case we have the following problems:

$\bu$ the spectrum of the Laplacian on infinite graphs usually has an absolutely continuous component;

$\bu$ due to the infiniteness of the graph there are infinitely many cycles of each length.

\medskip

Let $\cG=(\cV,\cE)$ be a $\G$-periodic graph with the quotient graph $\cG_*=\cG/\G=(\cV_*,\cE_*)$ (see Subsection \ref{ss1.1}). Atiyah \cite{A76} introduced the von Neumann trace of a $\G$-periodic bounded operator $T$ on $\ell^2(\cV)$ by
\[\lb{GTr}
\Tr_\G T=\sum_{v\in\cV_*}(T\d_v)_v,
\]
where $\d_v$ is an orthonormal basis of $\ell^2(\cV)$. In order to introduce and study the Ihara zeta function for periodic graphs, in \cite{GIL08,GIL09,LPS19} the $\G$-trace of the $n$-th power of Laplacians on periodic graphs was expressed  in terms of numbers of cycles of the \emph{periodic graph} of length $n$ starting at vertices from the fundamental domain of the periodic graph. According to the Floquet theory a $\G$-periodic operator $T$ can be decomposed into a direct integral of fiber operators $T(k)$, $k\in\T^d=\R^d/\Z^d$, acting on the finite quotient graph $\cG_*$, and the following identity holds true
\[\lb{DIGT}
\Tr_\G T=\frac1{(2\pi)^d}\int_{\T^d}\Tr T(k)\,dk,
\]
see \cite{A76}, \cite{KoS03}. The parameter $k$ is called the \emph{quasimomentum}. Trace formulas for the fiber adjacency operator on the quotient graph were discussed in \cite{AS87,S86}, where the traces of $n$-th power of the fiber adjacency operator were expressed in terms of the unitary characters of the quotient graph cycles with length $n$. In those papers trace formulas were used to describe the Ihara zeta function and to obtain the determinant identity for the fiber adjacency operator.

We do not know results about the trace formulas for the Schr\"odinger operators on periodic graphs. We study the traces for Laplace and Schr\"odinger operators on periodic graphs using the approach based on the decomposition of these operators into the direct integral of fiber operators acting on finite quotient graphs. Our goal is to obtain some identities connecting the band functions (which are the functions of the quasimomentum) of the Schr\"odinger operators on periodic graphs with geometry of the graphs. Therefore we introduce a modified quotient graph and reduce the fiber Schr\"odinger operators to the weighted adjacency operators on this modified graph. We also introduce the notion of \emph{cycle index} and divide all cycles of the quotient graph into disjoint subsets $\cC_n^m$ of cycles having the length $n$ and the index $m\in\Z^d$. Cycles of the quotient graph with zero index correspond to cycles in~the periodic graph. Cycles with non-zero indices are obtained by factorization of paths in the $\G$-periodic graph $\cG$ connecting some $\G$-equivalent vertices. Firstly, using a specific representation of the fiber operator \cite{KS14}, we express the traces of the
$n$-th power of the fiber operator as finite Fourier series of the quasimomentum. The coefficients of the Fourier series are given in terms of the potentials and cycles in the quotient graph from the set $\cC_n^m$. Here we essentially use that the cycle indices form the lattice $\Z^d$. Secondly, integrating these traces of fiber operators over the Brillouin zone, we obtain regularized traces which, due to \er{DIGT}, coincide with the $\G$-traces of the Laplace and Schr\"odinger operators on periodic graphs. The regularized traces are expressed in terms of cycles of the \emph{finite quotient graph} with zero index and the electric potential.

\medskip

The spectrum of the Schr\"odinger operator $H$ with a periodic potential $V$ on periodic graphs consists of a finite number of bands $\s_j$, $j=1,\ldots,\n$. We used the Fourier series for the traces of the $n$-th power of the fiber Schr\"odinger operators to prove the following new results:

$\bu$ the \emph{lower} estimates of the total bandwidth $\gS(H)= \sum_{j=1}^{\n}|\s_j|$ for the Schr\"odinger operator $H$ in terms of geometric parameters of the graph and the potential \cite{KS21a}:
$$
\gS(H)\geq\frac{2d}{v_*^{\n-1}}\,,
$$
where
$$
v_*=\vk_++\diam(V-\vk),\qqq \diam w=\max_{x}w_x-\min_{x}w_x,
$$
$\vk_+$ is the maximum vertex degree, and $\vk$ is the degree potential.

$\bu$ the spectrum of the magnetic Schr\"odinger operator with
each magnetic potential $t\a$, where $t$ is a coupling constant, has
an absolutely continuous component for all except finitely many $t$
from any bounded interval, see \cite{KS21b}.

Note that in our earlier papers \cite{KS14,KS17,KS18,KS20} we obtained only some \emph{upper} estimates of the total bandwidth for the Schr\"odinger operators. It is usually more difficult to obtain lower estimates.

\subsection{Schr\"odinger operators on periodic graphs.} \lb{ss1.1}
Let $\cG=(\cV,\cE)$ be a connected infinite graph, possibly  having
loops and multiple edges and embedded into the space $\R^d$. Here
$\cV$  is the set of its vertices and $\cE$ is the set of its
unoriented edges. Considering each edge in $\cE$ to have two
orientations, we introduce the set $\cA$ of all oriented edges. An
edge starting at a vertex $x$ and ending at a vertex $y$ from $\cV$
will be denoted as the ordered pair $(x,y)\in\cA$ and is said to be
\emph{incident} to the vertices. Let $\ul\be=(y,x)$ be the inverse
edge of $\be=(x,y)\in\cA$. Vertices $x,y\in\cV$ will be called
\emph{adjacent} and denoted by $x\sim y$, if $(x,y)\in \cA$. We
define the degree ${\vk}_x$ of the vertex $x\in\cV$ as the number of
all edges from $\cA$, starting at $x$.

Let $\G$ be a lattice of rank $d$ in $\R^d$ with a basis $\{a_1,\ldots,a_d\}$, i.e.,
$$
\G=\Big\{a : a=\sum_{s=1}^dn_sa_s, \; (n_s)_{s=1}^d\in\Z^d\Big\},
$$
and let
\[\lb{fuce}
\Omega=\Big\{\bx\in\R^d : \bx=\sum_{s=1}^d\bx_sa_s, \; (\bx_s)_{s=1}^d\in
[0,1)^d\Big\}
\]
be the \emph{fundamental cell} of the lattice $\G$. We define the equivalence relation on $\R^d$:
$$
\bx\equiv \by \; (\hspace{-4mm}\mod \G) \qq\Leftrightarrow\qq \bx-\by\in\G \qqq
\forall\, \bx,\by\in\R^d.
$$

We consider \emph{locally finite $\G$-periodic graphs} $\cG$, i.e.,
graphs satisfying the following conditions:
\begin{itemize}
  \item[1)] $\cG=\cG+a$ for any $a\in\G$;
  \item[2)] the quotient graph  $\cG_*=\cG/\G$ is finite.
\end{itemize}
The basis $a_1,\ldots,a_d$ of the lattice $\G$ is called the {\it
periods}  of $\cG$. We also call the quotient graph $\cG_*=\cG/\G$
the \emph{fundamental graph} of the periodic graph $\cG$. The
fundamental graph $\cG_*$ is a graph on the $d$-dimensional torus
$\R^d/\G$. The graph $\cG_*=(\cV_*,\cE_*)$ has the vertex set
$\cV_*=\cV/\G$, the set $\cE_*=\cE/\G$ of unoriented edges and the
set $\cA_*=\cA/\G$ of oriented edges which are finite.

Let $\ell^2(\cV)$ be the Hilbert space of all square summable
functions  $f:\cV\to \C$ equipped with the norm
$$
\|f\|^2_{\ell^2(\cV)}=\sum_{x\in\cV}|f_x|^2<\infty.
$$

We consider the Schr\"odinger operator $H$ acting on the Hilbert
space $\ell^2(\cV)$ and given by
\[\lb{Sh}
H=-\D+V,\qqq \D=\vk-A,
\]
where $\D$ is the combinatorial Laplacian, and $A$ is the
\emph{adjacency} operator having the form
\[
\lb{ALO}
(A f)_x=\sum_{(x,y)\in\cA}f_y, \qqq f\in\ell^2(\cV), \qqq x\in\cV,
\]
$V$ is a real $\G$-periodic potential, i.e., it satisfies for
all $(x,a)\in\cV\ts\G$:
\[
\lb{Pot} (Vf)_x=V_xf_x, \qqq V_{x+a}=V_x,
\]
and $\vk$ is a degree potential given by
\[\lb{DMO}
(\vk f)_x=\vk_xf_x, \qqq f\in\ell^2(\cV), \qqq x\in\cV,
\]
$\vk_x$ is the degree of the vertex $x$. It is known that the
Schr\"odinger operator $H$  is  a bounded self-adjoint operator on
$\ell^2(\cV)$ (see, e.g., \cite{SS92}) and the spectra of $A$ and
$\D$ satisfy
\[\lb{spAD}
\s(A)\subseteq[-\vk_+,\vk_+],\qqq
0\in\s(\D)\subseteq[0,2\vk_+],\qqq\where \qqq
\vk_+:=\max\limits_{x\in\cV_*}\vk_x.
\]

\begin{remarks}
1) If $\cG$ is a regular graph of degree $\vk_+$,
i.e., all vertices of $\cG$ have the same degree $\vk_+$, then the
Laplacian $\D$ has the the form
\[\lb{suD}
\D=\vk_+I-A,
\]
where $I$ is the identity operator, and $A$ is the adjacency
operator  given by \er{ALO}. Thus, the operators $-\D$ and $A$ on a
regular graph differ only by a shift.

2) When we consider the combinatorial Laplacian, without loss of
generality we may assume that there are no loops in the graph $\cG$.
Indeed, let $\wt\cG=(\cV,\wt\cA\,)$ be a graph with loops, and let
$\cG=(\cV,\cA)$ be a graph obtained from $\wt\cG$ by deleting all
its loops. Then the combinatorial Laplacians $\D_\cG$ and
$\D_{\wt\cG}$ on the graphs $\cG$ and $\wt\cG$, respectively, act on
the same space $\ell^2(\cV)$ and, due to \er{Sh} and \er{ALO},
$$
(\D_{\wt\cG}f)_x=\sum_{(x,y)\in\wt\cA}(f_x-f_y)=\sum_{(x,y)\in\cA}(f_x-f_y)=(\D_\cG
f)_x, \qqq \forall\, f\in\ell^2(\cV),\qqq \forall\, x\in\cV.
$$
Note that for the normalized Laplacian defined by \er{DNLA} this is
not true anymore.
\end{remarks}

\subsection{Edge indices} \lb{Sedin} Here we define an {\it edge index}.
It was introduced in \cite{KS14} in order to express fiber Laplacians and Schr\"odinger operators on the fundamental graph in terms of indices
(see \er{fado}) and estimate band widths.

For each $\bx\in\R^d$ we introduce the vector $\bx_\A\in\R^d$ by
\[\lb{cola}
\bx_\A=(\bx_1,\ldots,\bx_d), \qqq \textrm{where} \qq \bx=
\textstyle\sum\limits_{s=1}^d\bx_sa_s,
\]
i.e., $\bx_\A$ is the coordinate vector of $\bx$ with respect to the
basis  $\A=\{a_1,\ldots,a_d\}$ of the lattice~$\G$.

For any vertex $x\in\cV$ of a $\G$-periodic graph $\cG$ the
following  unique representation holds true:
\[\lb{Dv}
x=x_0+[x], \qq \textrm{where}\qq x_0\in\Omega,\qquad [x]\in\G,
\]
and $\Omega$ is a fundamental cell of the lattice $\G$ defined by
\er{fuce}. In other words, each vertex $x$ can be obtained from a
vertex $x_0\in \Omega$  by a shift by a vector $[x]\in\G$. For any
oriented edge $\be=(x,y)\in\cA$ we define the \emph{edge index}
$\t(\be)$ as the vector of the lattice $\Z^d$ given by
\[
\lb{in}
\t(\be)=[y]_\A-[x]_\A\in\Z^d,
\]
where $[x]\in\G$ is defined by \er{Dv} and the vector
$[x]_\A\in\Z^d$  is given by \er{cola}.

On the set $\cA$ of all oriented edges of the $\G$-periodic graph
$\cG$  we define the surjection
\[\lb{sur}
\gf:\cA\rightarrow\cA_*=\cA/\G,
\]
which maps each $\be\in\cA$ to its equivalence class
$\be_*=\gf(\be)$  which is an oriented edge of the fundamental graph
$\cG_*$. For any oriented edge $\be_*\in\cA_*$  we define the
\emph{edge index}  $\t(\be_*)\in\Z^d$ by
\[
\lb{dco} \t(\be_*)=\t(\be) \qq \textrm{ for some $\be\in\cA$ \; such
that }   \; \be_*=\gf(\be), \qqq \be_*\in\cA_*,
\]
where $\gf$ is defined by \er{sur}. In other words, edge indices of
the fundamental graph $\cG_*$ are induced by edge indices of the
periodic graph~$\cG$. The edge index $\t(\be_*)$ is uniquely
determined by \er{dco}, since
$$
\t(\be+a)=\t(\be),\qqq \forall\, (\be,a)\in\cA \ts\G.
$$
From the definition of the edge indices it follows that
\[\lb{inin}
\t(\ul\be\,)=-\t(\be), \qqq \forall\,\be\in\cA_*.
\]

\subsection{Direct integral decomposition.} We introduce the Hilbert space
\[\lb{Hisp}
\mH=L^2\Big(\T^{d},{dk\/(2\pi)^d}\,,\ell^2(\cV_*)\Big)
=\int_{\T^{d}}^{\os}\ell^2(\cV_*)\,{dk \/(2\pi)^d}\,, \qqq
\T^d=\R^d/(2\pi\Z)^d,
\]
i.e., a constant fiber direct integral, equipped with the norm
$$
\|g\|^2_{\mH}=\int_{\T^d}\|g(k)\|_{\ell^2(\cV_*)}^2\frac{dk}{(2\pi)^d}\,,
$$
where the function $g(k)\in\ell^2(\cV_*)$ for almost all $k\in\T^d$.
We recall Theorem 1.1 from \cite{KS14}.

\begin{theorem}\label{TFDC}
The Schr\"odinger operator $H=-\D+V$ on
$\ell^2(\cV)$ has  the following decomposition into a constant fiber
direct integral
\[
\lb{raz}
\begin{aligned}
& UH U^{-1}=\int^\oplus_{\T^d}H(k){dk\/(2\pi)^d}\,,
\end{aligned}
\]
where $U:\ell^2(\cV)\to\mH$ is some unitary operator (the Gelfand
transform), and the fiber Schr\"odinger operator $H(k)$ on
$\ell^2(\cV_*)$ is given by
\[\label{Hvt'}
H(k)=-\D(k)+V, \qqq \D(k)=\vk-A(k), \qqq \forall\,k\in\T^d.
\]
Here $V$ and $\vk$ are the electric and degree potentials  on
$\ell^2(\cV_*)$; $\D(k)$ is the fiber Laplacian, and $A(k)$ is the
fiber adjacency operator having the form
\[
\label{fado}
\big(A(k)f\big)_x=\sum_{\be=(x,y)\in\cA_*} e^{i\lan\t(\be),\,k\ran}f_y,
 \qqq f\in\ell^2(\cV_*),\qqq x\in \cV_*,
\]
where $\t(\be)$ is the index of the edge $\be\in\cA_*$ defined by \er{in}, \er{dco}.
\end{theorem}

\begin{remark}
The fiber operator $A(k)$ is expressed in terms of only edge indices
$\t(\be)$, $\be\in\cA_*$.
\end{remark}

\subsection{Spectrum of Schr\"odinger operators.} We briefly describe
the spectrum of the Schr\"o\-din\-ger operator $H$ (for more details
see, e.g., \cite{HS99} or \cite{KS14}). Let $\#M$ denote the number
of elements in a set $M$. Each fiber operator $H(k)$, $k\in\T^{d}$,
acts on the space $\ell^2(\cV_*)=\C^\n$, $\n=\#\cV_*$, and has $\n$
real eigenvalues $\l_j(k)$, $j=1,\ldots,\n$, labeled in
non-decreasing order by
\[
\label{eq.3H} \l_{1}(k)\leq\l_{2}(k)\leq\ldots\leq\l_{\nu}(k), \qqq
\forall\,k\in\T^{d},
\]
counting multiplicities.  Each $\l_j(\cdot)$ is  a real and
piecewise analytic function  on the torus $\T^{d}$ and creates the
\emph{spectral band} $\s_j(H)$ given by
\[\lb{ban.1H}
\s_j(H)=[\l_j^-,\l_j^+]=\l_j(\T^{d}),\qqq j\in\N_\n, \qqq \N_\n=\{1,\ldots,\n\}.
\]
Some of $\l_j(\cdot)$ may be constant, i.e.,
$\l_j(\cdot)=\L_j=\const$,  on some subset $\cB$ of $\T^d$ of
positive Lebesgue measure. In this case the Schr\"odinger operator
$H$ on $\cG$ has the eigenvalue $\L_j$ of infinite multiplicity. We
call $\{\L_j\}$ a \emph{flat band}. Thus, the spectrum of the
Schr\"odinger operator $H$ on the periodic graph $\cG$ has the form
\[\lb{specH}
\s(H)=\bigcup_{k\in\T^d}\s\big(H(k)\big)=
\bigcup_{j=1}^{\nu}\s_j(H)=\s_{ac}(H)\cup \s_{fb}(H),
\]
where $\s_{ac}(H)$ is the absolutely continuous spectrum, which is a
union of non-degenerate bands from \er{ban.1H}, and $\s_{fb}(H)$ is
the set of all flat bands (eigenvalues of infinite multiplicity).

\medskip

The paper is organized as follows. In Section \ref{Sec2} we formulate our main results:

$\bu$ trace formulas for the adjacency operators on periodic graphs (Theorem \ref{TFao});

$\bu$ trace formulas for the Schr\"odinger operators with periodic potentials on periodic graphs  (Theorems \ref{TPG} and \ref{TPG2});

$\bu$ trace formulas for the heat kernel and the resolvent of
the Schr\"odinger operators (Corollaries \ref{CHC} and \ref{CHC1}).

Section \ref{Sec3} is devoted to the proof of the main results. In Section \ref{Sec4} we consider some examples and apply the obtained results to calculate the traces for the adjacency and Schr\"odinger operators on the square lattice and on the Kagome lattice. Section \ref{Sec5} is devoted to the trace formulas for the normalized Laplacians on periodic graphs.

\section{Main results}
\setcounter{equation}{0}
\lb{Sec2}

\subsection{Cycle indices and cycle sets}\lb{SoC} In order to
formulate trace formulas we need some notation. A \emph{path} $\bp$
in a graph $\cG=(\cV,\cA)$ is a sequence of consecutive edges
\[\lb{depa}
\bp=(\be_1,\be_2,\ldots,\be_n), \qqq \where \qq
\be_s=(x_{s-1},x_s)\in\cA,\qq s\in \N_n,
\]
for some vertices $x_0,x_1,\ldots,x_n\in\cV$. The vertices $x_0$ and
$x_n$ are called the \emph{initial} and \emph{terminal} vertices of
the path $\bp$, respectively. If $x_0=x_n$, then the path $\bp$ is
called a \emph{cycle}. The number $n$ of edges in a cycle $\bc$ is
called the \emph{length} of $\bc$ and is denoted by $|\bc|$, i.e.,
$|\bc|=n$. The \emph{reverse} of the path $\bp$ given by \er{depa}
is the path $\ul\bp=(\ul\be_n,\ldots,\ul\be_1)$.

\begin{remark} A path $\bp$ is uniquely defined by the sequence of
its oriented edges $(\be_1,\be_2,\ldots,\be_n)$. The sequence of
its vertices $(x_0,x_1,\ldots,x_n)$ does not uniquely define~$\bp$,
since multiple edges are allowed in the graph $\cG$.
\end{remark}

We extend the notion of an \emph{edge index} defined in Subsection
\ref{Sedin} to a \emph{cycle index}. Let $\cC$ be the set of all
cycles of the fundamental graph $\cG_*$.  For any cycle $\bc\in\cC$
we define the \emph{cycle index} $\t(\bc)\in\Z^d$ by
\[\lb{cyin}
\t(\bc)=\sum\limits_{\be\in\bc}\t(\be),  \qqq \bc\in\cC.
\]
From this definition and the identity \er{inin}, it follows that
\[\lb{ininc}
\t(\ul\bc\,)=-\t(\bc), \qqq \forall\,\bc\in\cC.
\]
For the set  $\cC$ of all cycles of the fundamental graph
 $\cG_*=(\cV_*,\cA_*)$  we define the following subsets of $\cC$
which will be used throughout this
paper:\\
$\bu$  $\cC_n$ is the set of all cycles of length $n$ in $\cG_*$ and
$\cN_n$ is their number:
\[\lb{cNn}
\cC_{n}=\{\bc\in\cC: |\bc|=n \}, \qqq \cN_{n}=\#\cC_{n}<\iy;
\]
$\bu$  $\cC^\mm$ is the set of all cycles with index $\mm\in\Z^d$ in
$\cG_*$:
\[\lb{cNm}
\cC^\mm=\{\bc\in\cC: \t(\bc)=\mm \};
\]
$\bu$  $\cC_n^\mm$ is the set of all cycles of length $n$ and with
index $\mm$ in $\cG_*$ and  $\cN_n^\mm$ is their number:
\[
\lb{cNnm} \cC_n^\mm=\cC_n\cap\cC^\mm=\{\bc\in\cC: |\bc| =n
\;\textrm{and}\;\t(\bc)=\mm\}, \qqq \cN_n^\mm=\#\cC_n^\mm<\iy.
\]
Note that the sets $\cC_n$ and $\cC_n^\mm$ are finite. Here and
below\\
$\bu$  $|\bc|$ is the length  of the cycle $\bc$;\\
 $\bu$  $\t(\bc)$ is the index of the cycle $\bc$ defined by \er{cyin}.

\begin{remarks}
1) Any cycle $\bc$ in the fundamental graph $\cG_*$ is obtained by
factorization of a path in the periodic graph $\cG$ connecting some
$\G$-equivalent vertices $x\in\cV$ and $x+a\in\cV$, $a\in\G$.
Furthermore, the index of the cycle $\bc$ is equal to
$\mm=(m_j)_{j=1}^d\in\Z^d$, where $a=m_1a_1+\ldots+m_da_d$. In
particular, $\t(\bc)=0$ if and only if the cycle $\bc$ in $\cG_*$
corresponds to a cycle in~$\cG$.

2) The sets $\cC,\cC_n,\cC^\mm,\cC_n^\mm$ include the corresponding
cycles  with \emph{back-tracking parts}, i.e., cycles
$(\be_1,\ldots,\be_n)$ for which $\be_{s+1}=\ul\be_s$ for some
$s\in\N_n$ ($\be_{n+1}$ is understood as $\be_1$).

3) Since there are no loops in the periodic graph $\cG$, the
fundamental  graph $\cG_*$ has no loops with zero index. Note that
$\cG_*$ may have loops with non-zero index. Since there are no loops
with zero index in $\cG_*$ and each backtracking contributes 2 in
the cycle length, all cycles of length 3 with zero index in $\cG_*$
are \emph{proper} cycles, i.e., cycles without backtracking.
\end{remarks}

We formulate some simple properties of the numbers $\cN_n$ and
$\cN_n^\mm$. A graph is called \emph{bipartite} if its vertex set is
divided into two disjoint sets (called \emph{parts} of the graph)
such that each edge connects vertices from distinct parts.

\begin{proposition}
\lb{spcN} i) The numbers $\cN_n=\#\cC_n$ and $\cN_n^\mm=\#\cC_n^\mm$
 satisfy
\[\lb{cN10}
\cN_1^0=0,\qqq \#\cA_*\le \cN_{2n}^0,
\]
\[\lb{Npme}
\cN_n\leq\n\vk^n_+,\qqq \cN_n^\mm=\cN_n^{-\mm},\qqq
\forall\,(n,\mm)\in\N\ts\Z^d,\qqq \n=\#\cV_*,
\]
\[\lb{cNe0}
\cN_n^\mm=0, \qqq \textrm{if}\qqq
\|\mm\|>n\t_+,\qqq\textrm{where}\qqq
\t_+=\max\limits_{\be\in\cA_*}\|\t(\be)\|,
\]
$\t(\be)$ is the edge index  defined by \er{in}, \er{dco}, and
$\|\cdot\|$ is the standard norm in $\R^d$; $\vk_+$ is defined in
\er{spAD}.

If there are no multiple edges in the periodic graph $\cG$, then
\[\lb{cN0g1}
\cN_2^0=\#\cA_*.
\]

ii) If the fundamental graph $\cG_*=(\cV_*,\cA_*)$ is bipartite,
then
\[\lb{cNBG}
\cN_n=0 \qq \textrm{and}\qq \cN_n^\mm=0 \qq \textrm{for odd \;
$n\in\N$ \; and all \; $\mm\in\Z^d$}.
\]
Moreover, the fundamental graph $\cG_*$ is bipartite iff $\cN_n=0$
for all odd $n\leq\n$.

iii) The periodic graph $\cG$ is bipartite iff
\[\lb{cN0BG}
\cN_n^0=0 \qq \textrm{for all odd}\qq n.
\]

iv) For any $s\in\N$ there exists a periodic graph such that
\[\lb{CNs}
\cN_n^0=0\qq \textrm{for all odd $n<2s$,\qq and} \qq
\cN_{2s+1}^0\neq0.
\]
\end{proposition}

\subsection{Trace formulas for adjacency operators}
First we discuss trace formulas for the fiber adjacency operator.
The eigenvalues of the fiber adjacency operator $A(k)$ will be
denoted  by $\l^o_j(k)$, $j\in\N_\n$. The spectral bands $\s_j(A)$,
$j\in\N_{\n}$, for the adjacency operator $A$ have the form
$$
\s_j(A)=[\l_j^{o-},\l_j^{o+}]=\l_j^o(\T^d).
$$
We describe trace formulas for the adjacency operator $A$ in terms
of the numbers $\cN_n$ and $\cN_n^m$.

\begin{theorem}
\lb{TFao} Let $A(k)$, $k\in\T^d$, be the fiber adjacency operator
defined by  \er{fado} on the fundamental graph $\cG_*=(\cV_*,\cA_*)$
with $\n=\#\cV_*$. Then for each $n\in\N$

i) The trace of $A^n(k)$ has the form
\[\lb{TrAo}
\Tr
A^n(k)=\sum_{j=1}^\n\big(\l_j^o(k)\big)^n=\sum_{\bc\in\cC_n}\cos\lan
\t(\bc),k\ran =\sum_{\mm\in\Z^d\atop\|\mm\|\leq
n\t_+}\cN_n^\mm\cos\lan\mm,k\ran.
\]
Moreover,
\[\lb{EsTAn}
\big|\Tr A^n(k)\big|\leq\Tr A^n(0)=\cN_n.
\]
Here $\t_+$ and $\vk_+$ are given in \er{cNe0} and \er{spAD},
respectively; $\|\cdot\|$ is the standard norm in $\R^d$.

ii) The trace of $A^n(k)$ satisfies
\[\lb{ITrA}
\frac1{(2\pi)^d}\int_{\T^d}\Tr A^n(k)dk=\cN_n^0\geq0.
\]

iii) The fundamental graph $\cG_*$ is bipartite iff
\[\lb{Kbfg}
\Tr A^n(0)=0\qqq \textrm{for all odd $n\leq\n$}.
\]

iv) The periodic graph $\cG$ is bipartite iff
\[\lb{Kbpg}
\int_{\T^d}\Tr A^n(k)dk=0\qqq \textrm{for all odd $n$}.
\]
Moreover, the condition \er{Kbpg} can not be reduced, i.e., for any $s\in\N$ there exists a non-bipartite periodic graph $\cG$ such that $\int_{\T^d}\Tr A^n(k)dk=0$ for all odd $n<s$.
\end{theorem}

\begin{remarks}
1) The identity $\Tr A^n(0)=\cN_n$, $n\in\N$, see \er{EsTAn}, is the
well-known trace formula for the adjacency operator $A(0)$ on the finite
graph $\cG_*$, see, e.g., \cite{Mn07}.

2) The fundamental graph of a bipartite periodic graph is not
necessary  bipartite, see Example \ref{ExBB}. But for any bipartite
periodic graph there exists a bipartite fundamental graph, see
Remark after Proof of Example \ref{ExSL}.
\end{remarks}

\subsection{Trace formulas for Schr\"odinger operators}
In order to determine trace formulas for the Schr\"odinger operators
$H=-\D+V$, we need some modification of the fundamental graph
$\cG_*$. At each vertex $x$ of $\cG_*$ we add a loop with weight
$v_x=V_x-\vk_x$. This modification allows us to replace the fiber
Schr\"odinger operator $H(k)$ acting on the fundamental graph
$\cG_*$ by a fiber weighted adjacency operator on the modified
fundamental graph, and to obtain trace formulas for $H(k)$ in the
form similar to \er{TrAo}.

More precisely, at each vertex $x$ of $\cG_*=(\cV_*,\cA_*)$ we add a
loop  $\be_x$ with index $\t(\be_x)=0$ and consider the modified
fundamental graph $\wt\cG_*=(\cV_*,\wt\cA_*)$, where
\[\lb{wtAs}
\wt\cA_*=\cA_*\cup\{\be_x\}_{x\in\cV_*}, \qqq \t(\be_x)=0.
\]

Let $\wt\cC$ be the set of all cycles in $\wt\cG_*$. For each cycle
$$
\bc=(\be_1,\ldots,\be_n)\in\wt\cC, \qqq  \be_j=(x_j,x_{j+1}),
\qq j\in \N_n, \qq x_{n+1}=x_1,
$$
we define the \emph{weight} $\o(\bc)$  and the sum $v(\bc)$ by
\[\lb{Wcy} \o(\bc)=\o(\be_1)\ldots \o(\be_n),\qqq
v(\bc)=v_{x_1}+\ldots+v_{x_n},
\]
where $\o(\be)$ has the form
\[\lb{webe}
\o(\be)=\left\{
\begin{array}{cl}
1,  & \qq \textrm{if} \qq  \be\in\cA_* \\[2pt]
v_x, & \qq \textrm{if} \qq \be=\be_x
\end{array}\right.,\qqq v_x=V_x-\vk_x,
\]
and $\vk_x$ is the degree of the vertex $x\in\cV_*$.

\begin{remark}
Note that
\[\lb{oc1C}
\o(\bc)=1\qqq \textrm{for each cycle}\qqq \bc\in\cC.
\]
\end{remark}

Let the cycle sets $\wt\cC_n$, $\wt\cC^\mm$, $\wt\cC_n^\mm$,
and the numbers $\wt\cN_n$ and $\wt\cN_n^\mm$
for the modified fundamental graph $\wt\cG_*$ be defined as the
corresponding cycle sets $\cC_n$, $\cC^\mm$, $\cC_n^\mm$,
and the numbers $\cN_n$ and $\cN_n^\mm$ for the fundamental
graph $\cG_*$ (see Subsection \ref{SoC}).

We define a functional $\cT_n(k)$ as the finite sum over cycles from $\wt\cC_n$ given by
\[
\lb{deTnk}
\cT_n(k)=\sum_{\bc\in\wt\cC_n}\o(\bc)e^{-i\lan\t(\bc),k\ran} =
\sum_{\bc\in\wt\cC_n}\o(\bc)\cos\lan\t(\bc),k\ran,
\]
where we have used that $\o(\bc)=\o(\ul\bc\,)$ and
$\t(\bc)=-\t(\ul\bc\,)$  for each $\bc\in\wt\cC_n$. The functional
$\cT_n(k)$ is a finite Fourier series of the quasimomentum $k\in
\T^d$, where the coefficients $\o(\bc)$ are polynomials of degree $\le
n$  with respect to the potential $V$. We sometimes write
$\cT_n(k,V),\o(\bc,V),\ldots\,$ instead of
$\cT_n(k),\o(\bc),\ldots\,$, when several potentials $V$ are dealt
with.

We formulate some properties of the functionals $\cT_n(k)$. Denote
by $\o_+$  the following number
\[\lb{obc+}
\o_+=\max_{\be\in\wt\cA_*}|\o(\be)|=\max\big\{1,\|v\|_{\ell^\iy(\cV_*)}\big\},
\qqq v=(v_x)_{x\in\cV_*}.
\]

\begin{proposition}\lb{Ppof}
Let the functional $\cT_n(k)$, $n\in\N$, $k\in\T^d$, be defined by \er{deTnk}. Then

i) $\cT_n(k)$  satisfies
\[
\lb{Wesgh}
\big|\cT_n(k)\big|\leq\o_+^n\,\wt\cN_n,
\]
where $\o_+$ is given in \er{obc+}.

ii) $\cT_n(k)$ has the following finite Fourier series
\[\lb{FsTrH}
\cT_n(k)=\sum_{\mm\in\Z^d\atop \|\mm\|\leq n\t_+}
\cT_{n,\mm}\cos\lan\mm,k\ran,\qqq
\cT_{n,\mm}=\sum_{\bc\in\wt\cC_n^\mm}\o(\bc),
\]
and
\[\lb{IVFu}
\frac1{(2\pi)^d}\int_{\T^d}\cT_n(k) dk=\cT_{n,0},\qqq
\cT_{n,0}=\sum_{\bc\in\wt\cC_n^0}\o(\bc),
\]
where $\t_+$ and $\o(\bc)$ are defined in \er{cNe0} and \er{Wcy},
respectively, and $\|\cdot\|$ is the standard norm in $\R^d$.
\end{proposition}

\medskip

We formulate trace formulas for the fiber Schr\"odinger operators.
Recall that \begin{itemize}
  \item $|\bc|$ is the length  of the cycle $\bc$;
  \item $\t(\bc)$ is the index of the cycle $\bc$ defined by \er{cyin};
  \item $\o(\bc)$ is the weight of the cycle $\bc$ defined in \er{Wcy};
  \item $v(\bc)$ is the sum of the potential values along the cycle $\bc$ defined in \er{Wcy}.
\end{itemize}

\begin{theorem}\lb{TPG}
Let $H(k)$, $k\in\T^d$, be the fiber Schr\"odinger operator defined
by \er{Hvt'} -- \er{fado} on the fundamental graph
$\cG_*=(\cV_*,\cA_*)$. Then for each $n\in\N$ the trace of $H^n(k)$ has the form
\[\lb{TrnH}
\Tr H^n(k)=\sum_{j=1}^\n\l_j^n(k)=\cT_n(k), \qqq \n=\#\cV_*,
\]
where $\cT_n(k)$ is defined by \er{deTnk}. In particular,
\[\lb{TrH123}
\begin{aligned}
&\Tr H(k)=\cT_1(k)=\sum_{x\in\cV_*}v_x+
\sum_{\bc\in\cC_1}\cos\lan\t(\bc),k\ran,\qqq v_x=V_x-\vk_x,\\
&\Tr H^2(k)=\cT_2(k)=\sum_{x\in\cV_*}v_x^2+
2\sum_{\bc\in\cC_1}v(\bc)\cos\lan\t(\bc),k\ran+
\sum_{\bc\in\cC_2}\cos\lan\t(\bc),k\ran.
\end{aligned}
\]
If there are no loops in the fundamental graph $\cG_*$, then
\[\lb{hnl}
\Tr H(k)=\cT_1(k)=\sum_{x\in\cV_*}v_x,
\]
and if, in addition, there are no multiple edges in $\cG_*$, then
\[\lb{nome}
\Tr H^2(k)=\cT_2(k)=\#\cA_*+\|v\|_*^2,
\]
where $\|v\|_*^2=\sum\limits_{x\in\cV_*}v_x^2$ and
$\#\cA_*=\sum\limits_{x\in\cV_*}\vk_x>0$.
\end{theorem}

We determine traces of the Schr\"odinger operators in the integral
forms.

\begin{theorem}\lb{TPG2}
Let $H(k)$, $k\in\T^d$, be the fiber Schr\"odinger operator defined
by \er{Hvt'} -- \er{fado} on the fundamental graph
$\cG_*=(\cV_*,\cA_*)$. Then for each $n\in\N$ the trace of $H^n(k)$
satisfies
\[\lb{ITrH}
\frac1{(2\pi)^d}\int_{\T^d}\Tr H^n(k)dk=\cT_{n,0},
\]
where $\cT_{n,0}$ is defined in \er{IVFu}. In particular,
\[\lb{TR13IH}
\begin{aligned}
&\frac1{(2\pi)^d}\int_{\T^d}\Tr H(k)dk=\cT_{1,0}=\sum_{x\in\cV_*}v_x,\qqq v_x=V_x-\vk_x,
\\
&\frac1{(2\pi)^d}\int_{\T^d}\Tr H^2(k)dk=\cT_{2,0}=\sum_{x\in\cV_*}
v_x^2+\cN_2^0.
\end{aligned}
\]
If there are no multiple edges in the periodic graph $\cG$, then
$\cN_2^0=\#\cA_*$ in \er{TR13IH}.
\end{theorem}

\begin{remarks}
1) The formulas \er{TrnH}, \er{deTnk} and \er{ITrH}, \er{IVFu} are
 \emph{trace formulas}, where the traces of the fiber operators are expressed in terms of
some geometric parameters of the graph (vertex degrees, cycle
indices and lengths) and the potential $V$.

2) The trace formulas for the Laplacian $-\D=A-\vk$ are given by the
same identities \er{TrnH} and \er{ITrH}, where
$\cT_n(k)=\cT_n(k,0)$ and $\cT_{n,0}=\cT_{n,0}(0)$.

3) The regularized traces \er{ITrH} of the operators $H^n(k)$ coincide with the $\G$-traces of $H^n$, see \er{GTr}, \er{DIGT}.
\end{remarks}

In the following statement we compare the traces of the fiber
Schr\"odinger  operators and the traces of the corresponding fiber
adjacency operators.

\begin{corollary}\lb{CCrFA}
Let $H(k)=A(k)-\vk+V$, $k\in\T^d$, be the fiber Schr\"odinger operator,
where $A(k)$ is the fiber adjacency operator given by \er{fado}. Then for each $n\in\N$
\[\label{TrHAn}
\Tr\big(H^n(k)-A^n(k)\big)
=\hspace{-3mm}\sum_{\bc\in\wt\cC_n\sm\cC_n}\hspace{-3mm}\o(\bc)\cos\lan\t(\bc),k\ran=
\hspace{-3mm}\sum_{\mm\in\Z^d \atop \|\mm\|\leq n\t_+}\hspace{-3mm}\gt_{n,\mm}\cos\lan\mm,k\ran,
\qq \gt_{n,\mm}=\hspace{-3mm}\sum_{\bc\in\wt\cC_n^\mm\sm\cC_n^\mm}\hspace{-3mm}\o(\bc),
\]
\[\label{TrHAnI}
\frac1{(2\pi)^d}\int_{\T^d}\Tr\big(H^n(k)-A^n(k)\big)dk
=\gt_{n,0}, \qqq \gt_{n,0}=
\sum_{\bc\in\wt\cC_n^0\sm\cC_n^0}\o(\bc),
\]
where $\t_+$ is defined in \er{cNe0}, and $\|\cdot\|$ is the
standard norm in $\R^d$.  In particular,
\[\label{TrHA1}
\begin{aligned}
&\Tr\big(H(k)-A(k)\big)=\sum_{x\in\cV_*}v_x,\qqq v_x=V_x-\vk_x,\\
&\Tr\big(H^2(k)-A^2(k)\big)=\|v\|_*^2+
2\sum_{\bc\in\cC_1}v(\bc)\cos\lan\t(\bc),k\ran, \qqq
\|v\|_*^2=\sum_{x\in\cV_*}v_x^2,
\end{aligned}
\]
\[\label{TrHA1I}
\begin{aligned}
&\frac{1}{(2\pi)^d}\int_{\T^d}\Tr\big(H(k)-A(k)\big)dk=\sum_{x\in\cV_*}v_x,
\\
&\frac1{(2\pi)^d}\int_{\T^d}\Tr\big(H^2(k)-A^2(k)\big)dk=
\|v\|_*^2\geq0.
\end{aligned}
\]
\end{corollary}

\begin{remarks}
1) From the second identity in \er{TrHA1I} it simply follows that
\[
\lb{HA} H=A \qqq \textrm{iff} \qqq
\int_{\T^d}\Tr\big(H^2(k)-A^2(k)\big)dk=0.
\]

2) The explicit form for the identities \er{TrnH} and \er{ITrH} as $n=3$ is given in Proposition~\ref{TrTO}.
\end{remarks}

\subsection{Traces of the heat kernel and the resolvent}
We present trace formulas for the heat kernel and the resolvent of
the fiber Schr\"odinger operators.

\begin{corollary} \lb{CHC} Let $H(k)$, $k\in\T^d$, be the fiber
Schr\"odinger operator defined by \er{Hvt'} -- \er{fado} on the
fundamental graph $\cG_*=(\cV_*,\cA_*)$. Then the trace of $e^{t
H(k)}$, $t\in\C$, satisfies:
\[\lb{hkfH}
\Tr e^{tH(k)}=\n+\sum_{n=1}^\iy\frac{t^n}{n!}\,\cT_n(k)=
\n+\sum_{\bc\in\wt\cC}\frac{\o(\bc)}{|\bc|!}\,
t^{|\bc|}\cos\lan\t(\bc),k\ran, \qqq \n=\#\cV_*,
\]
\[\lb{THFS}
\Tr e^{tH(k)}=\n+\sum\limits_{\mm\in\Z^d}\gh_\mm(t)\cos\lan\mm,k\ran,
\qqq \gh_\mm(t)=\sum_{\bc\in\wt\cC^\mm}\frac{\o(\bc)}{|\bc|!}\,
t^{|\bc|},
\]
\[\lb{hkfHI}
\frac1{(2\pi)^d}\int_{\T^d}\Tr e^{tH(k)}dk=\n+\sum_{n=1}^\iy\frac{t^n}{n!}\,\cT_{n,0}=
\n+\sum_{\bc\in\wt\cC^0}\frac{\o(\bc)}{|\bc|!}\,
t^{|\bc|},
\]
where $\cT_n(k)$ and $\cT_{n,0}$ are given in \er{deTnk} and \er{IVFu},
respectively. The series in \er{hkfH} -- \er{hkfHI} converge absolutely for all $t\in\C$.
\end{corollary}

\begin{remarks}
1) From the expansion \er{hkfH} we see that the functionals $\cT_n(k)$ defined by \er{deTnk} are the \emph{heat kernel coefficients}.

2) The first heat kernel coefficients $\cT_n(k)$ and $\cT_{n,0}$, $n=1,2$,
 are given in \er{TrH123} and \er{TR13IH}, respectively. For $n=3$ see Proposition \ref{TrTO}.

3) The heat kernel for the Laplacian on regular graphs was considered in \cite{CJK15}.
\end{remarks}

\begin{corollary} \lb{CHC1} Let $H(k)$, $k\in\T^d$, be the fiber Schr\"odinger
 operator defined by \er{Hvt'} -- \er{fado} on the fundamental graph
$\cG_*=(\cV_*,\cA_*)$ with $\n=\#\cV_*$. Then the trace of the resolvent of $H(k)$ has
the following expansions:
\begin{multline}\label{reexH0}
\Tr\big(H(k)-\l I\big)^{-1}=
-\frac\n\l-\sum_{n=1}^\iy\frac{\cT_n(k)}{\l^{n+1}}
=-\frac\n\l-\sum_{\bc\in\wt\cC}
\frac{\o(\bc)}{\l^{|\bc|+1}}\,e^{-i\lan\t(\bc),k\ran}
\\ =-\frac\n\l-\sum_{\bc\in\wt\cC}\frac{\o(\bc)}{\l^{|\bc|+1}}\,\cos\lan\t(\bc),k\ran,
\end{multline}
\[\lb{reexHF}
\Tr\big(H(k)-\l I\big)^{-1}=-\frac\n\l-\sum_{\mm\in\Z^d}\cR_\mm(\l)\cos\lan\mm,k\ran,
\qq\textrm{where}\qq
\cR_\mm(\l)=\sum_{\bc\in\wt\cC^\mm}
\frac{\o(\bc)}{\l^{|\bc|+1}}\,,
\]
\[\lb{TrRI}
\frac1{(2\pi)^d}\int_{\T^d}\Tr \big(H(k)-\l I\big)^{-1}dk=
-\frac\n\l-\sum_{n=1}^\iy\frac{\cT_{n,0}}{\l^{n+1}}
=-\frac\n\l-\sum_{\bc\in\wt\cC^0}\frac{\o(\bc)}{\l^{|\bc|+1}},
\]
where $\cT_n(k)$ and $\cT_{n,0}$ are given in \er{deTnk} and
\er{IVFu}, respectively. The series in \er{reexH0}, \er{reexHF}
converge absolutely for $|\l|>\|H(k)\|$, and the series in \er{TrRI} converges absolutely for $|\l|>\|H\|$.
\end{corollary}

\begin{remark}
The first coefficients $\cT_n(k)$ and $\cT_{n,0}$, $n=1,2$, of the expansions \er{reexH0} and \er{TrRI} are given in \er{TrH123} and \er{TR13IH}, respectively. For $n=3$ see Proposition \ref{TrTO}.
\end{remark}

\subsection{Example} We apply the obtained results to the Schr\"odinger operator with
periodic potentials on the square lattice.

\begin{figure}[h]
\centering
\unitlength 1.0mm 
\linethickness{0.4pt}
\ifx\plotpoint\undefined\newsavebox{\plotpoint}\fi 
\begin{picture}(130,43)(0,0)
\put(23.5,24){$\Omega$}
\put(5,10){\line(1,0){40.00}}
\put(5,20){\line(1,0){40.00}}
\put(5,30){\line(1,0){40.00}}
\put(5,40){\line(1,0){40.00}}
\put(10,5){\line(0,1){38.00}}
\put(20,5){\line(0,1){38.00}}
\put(30,5){\line(0,1){38.00}}
\put(40,5){\line(0,1){38.00}}

\bezier{30}(10.5,10)(10.5,20)(10.5,30)
\bezier{30}(11.0,10)(11.0,20)(11.0,30)
\bezier{30}(11.5,10)(11.5,20)(11.5,30)
\bezier{30}(12.0,10)(12.0,20)(12.0,30)
\bezier{30}(12.5,10)(12.5,20)(12.5,30)
\bezier{30}(13.0,10)(13.0,20)(13.0,30)
\bezier{30}(13.5,10)(13.5,20)(13.5,30)
\bezier{30}(14.0,10)(14.0,20)(14.0,30)
\bezier{30}(14.5,10)(14.5,20)(14.5,30)
\bezier{30}(15.0,10)(15.0,20)(15.0,30)
\bezier{30}(15.5,10)(15.5,20)(15.5,30)
\bezier{30}(16.0,10)(16.0,20)(16.0,30)
\bezier{30}(16.5,10)(16.5,20)(16.5,30)
\bezier{30}(17.0,10)(17.0,20)(17.0,30)
\bezier{30}(17.5,10)(17.5,20)(17.5,30)
\bezier{30}(18.0,10)(18.0,20)(18.0,30)
\bezier{30}(18.5,10)(18.5,20)(18.5,30)
\bezier{30}(19.0,10)(19.0,20)(19.0,30)
\bezier{30}(19.5,10)(19.5,20)(19.5,30)
\bezier{30}(20.0,10)(20.0,20)(20.0,30)
\bezier{30}(20.5,10)(20.5,20)(20.5,30)
\bezier{30}(21.0,10)(21.0,20)(21.0,30)
\bezier{30}(21.5,10)(21.5,20)(21.5,30)
\bezier{30}(22.0,10)(22.0,20)(22.0,30)
\bezier{30}(22.5,10)(22.5,20)(22.5,30)
\bezier{30}(23.0,10)(23.0,20)(23.0,30)
\bezier{30}(23.5,10)(23.5,20)(23.5,30)
\bezier{30}(24.0,10)(24.0,20)(24.0,30)
\bezier{30}(24.5,10)(24.5,20)(24.5,30)
\bezier{30}(25.0,10)(25.0,20)(25.0,30)
\bezier{30}(25.5,10)(25.5,20)(25.5,30)
\bezier{30}(26.0,10)(26.0,20)(26.0,30)
\bezier{30}(26.5,10)(26.5,20)(26.5,30)
\bezier{30}(27.0,10)(27.0,20)(27.0,30)
\bezier{30}(27.5,10)(27.5,20)(27.5,30)
\bezier{30}(28.0,10)(28.0,20)(28.0,30)
\bezier{30}(28.5,10)(28.5,20)(28.5,30)
\bezier{30}(29.0,10)(29.0,20)(29.0,30)
\bezier{30}(29.5,10)(29.5,20)(29.5,30)

\put(10,10){\circle*{1}}
\put(20,10){\circle*{1}}
\put(30,10){\circle{1}}
\put(40,10){\circle{1}}

\put(10,20){\circle*{1}}
\put(20,20){\circle*{1}}
\put(30,20){\circle{1}}
\put(40,20){\circle{1}}

\put(10,30){\circle{1}}
\put(20,30){\circle{1}}
\put(30,30){\circle{1}}
\put(40,30){\circle{1}}

\put(10,40){\circle{1}}
\put(20,40){\circle{1}}
\put(30,40){\circle{1}}
\put(40,40){\circle{1}}

\put(10,10){\vector(1,0){20.00}}
\put(10,10){\vector(0,1){20.00}}
\put(7.0,8.0){$\scriptstyle x_1$}
\put(17.0,8.0){$\scriptstyle x_3$}
\put(30.5,8.0){$\scriptstyle x_1+a_1$}
\put(30.5,18.0){$\scriptstyle x_4+a_1$}
\put(26.5,31.0){$\scriptstyle x_1+a_1+a_2$}
\put(7.0,18.0){$\scriptstyle x_4$}
\put(17.0,18.0){$\scriptstyle x_2$}
\put(7.0,31.0){$\scriptstyle x_1+a_2$}
\put(17.0,31.0){$\scriptstyle x_3+a_2$}
\put(5,35){$\dL^2$}
\put(25,8){$\scriptstyle a_1$}
\put(6.5,26){$\scriptstyle a_2$}

\put(0,5){\emph{a})}
\put(80.5,26){$\Omega$}
\put(87,38){$\dL_*^2$}
\put(65,10){\line(1,0){20.00}}
\put(65,10){\line(0,1){20.00}}
\put(65,20){\line(1,0){20.00}}
\put(75,10){\line(0,1){20.00}}

\bezier{30}(65.5,10)(65.5,20)(65.5,30)
\bezier{30}(66.0,10)(66.0,20)(66.0,30)
\bezier{30}(66.5,10)(66.5,20)(66.5,30)
\bezier{30}(67.0,10)(67.0,20)(67.0,30)
\bezier{30}(67.5,10)(67.5,20)(67.5,30)
\bezier{30}(68.0,10)(68.0,20)(68.0,30)
\bezier{30}(68.5,10)(68.5,20)(68.5,30)
\bezier{30}(69.0,10)(69.0,20)(69.0,30)
\bezier{30}(69.5,10)(69.5,20)(69.5,30)
\bezier{30}(70.0,10)(70.0,20)(70.0,30)
\bezier{30}(70.5,10)(70.5,20)(70.5,30)
\bezier{30}(71.0,10)(71.0,20)(71.0,30)
\bezier{30}(71.5,10)(71.5,20)(71.5,30)
\bezier{30}(72.0,10)(72.0,20)(72.0,30)
\bezier{30}(72.5,10)(72.5,20)(72.5,30)
\bezier{30}(73.0,10)(73.0,20)(73.0,30)
\bezier{30}(73.5,10)(73.5,20)(73.5,30)
\bezier{30}(74.0,10)(74.0,20)(74.0,30)
\bezier{30}(74.5,10)(74.5,20)(74.5,30)
\bezier{30}(75.0,10)(75.0,20)(75.0,30)
\bezier{30}(75.5,10)(75.5,20)(75.5,30)
\bezier{30}(76.0,10)(76.0,20)(76.0,30)
\bezier{30}(76.5,10)(76.5,20)(76.5,30)
\bezier{30}(77.0,10)(77.0,20)(77.0,30)
\bezier{30}(77.5,10)(77.5,20)(77.5,30)
\bezier{30}(78.0,10)(78.0,20)(78.0,30)
\bezier{30}(78.5,10)(78.5,20)(78.5,30)
\bezier{30}(79.0,10)(79.0,20)(79.0,30)
\bezier{30}(79.5,10)(79.5,20)(79.5,30)
\bezier{30}(80.0,10)(80.0,20)(80.0,30)
\bezier{30}(80.5,10)(80.5,20)(80.5,30)
\bezier{30}(81.0,10)(81.0,20)(81.0,30)
\bezier{30}(81.5,10)(81.5,20)(81.5,30)
\bezier{30}(82.0,10)(82.0,20)(82.0,30)
\bezier{30}(82.5,10)(82.5,20)(82.5,30)
\bezier{30}(83.0,10)(83.0,20)(83.0,30)
\bezier{30}(83.5,10)(83.5,20)(83.5,30)
\bezier{30}(84.0,10)(84.0,20)(84.0,30)
\bezier{30}(84.5,10)(84.5,20)(84.5,30)
\bezier{30}(85.0,10)(85.0,20)(85.0,30)

\put(65,10){\circle*{1}}
\put(65,20){\circle*{1}}
\put(65,30){\circle{1}}
\put(75,10){\circle*{1}}
\put(75,20){\circle*{1}}
\put(75,30){\circle{1}}
\put(85,10){\circle{1}}
\put(85,20){\circle{1}}
\put(85,30){\circle{1}}

\put(62,8){$\scriptstyle x_1$}
\put(85.5,8){$\scriptstyle x_1$}
\put(62,31.0){$\scriptstyle x_1$}
\put(85.0,31.0){$\scriptstyle x_1$}

\put(65,10){\vector(1,0){20.00}}
\put(65,10){\vector(0,1){20.00}}

\put(73.7,8){$\scriptstyle x_3$}
\put(62.0,18){$\scriptstyle x_4$}
\put(72,18){$\scriptstyle x_2$}
\put(74,31.0){$\scriptstyle x_3$}
\put(85.3,18){$\scriptstyle x_4$}
\put(79,8){$\scriptstyle a_1$}
\put(61.5,25){$\scriptstyle a_2$}

\put(69,21.0){$\scriptstyle \be_2$}
\put(79,21.0){$\scriptstyle \be_6$}
\put(69,11.0){$\scriptstyle \be_4$}
\put(79,11.0){$\scriptstyle \be_8$}
\put(65.5,25){$\scriptstyle \be_5$}
\put(75.5,25){$\scriptstyle \be_7$}
\put(65.5,15){$\scriptstyle \be_1$}
\put(75.5,15){$\scriptstyle \be_3$}

\multiput(65,30)(4,0){5}{\line(1,0){2}}
\multiput(85,10)(0,4){5}{\line(0,1){2}}

\put(56,5){\emph{b})}


\put(100,10){\line(1,0){20.00}}
\put(100,10){\line(0,1){20.00}}
\put(100,30){\line(1,0){20.00}}
\put(120,10){\line(0,1){20.00}}
\put(100,10){\circle*{1}}
\put(100,30){\circle*{1}}
\put(120,10){\circle*{1}}
\put(120,30){\circle*{1}}

\put(96.5,7.5){$x_1$}
\put(120.5,7.5){$x_3$}
\put(96.5,31.5){$x_4$}
\put(120.0,31.0){$x_2$}

\bezier{200}(100,10)(110,2)(120,10)
\bezier{200}(100,30)(110,38)(120,30)
\bezier{200}(100,10)(92,20)(100,30)
\bezier{200}(120,10)(128,20)(120,30)
\put(110,6){\vector(1,0){1.00}}
\put(110,34){\vector(-1,0){1.00}}
\put(96,20){\vector(0,-1){1.00}}
\put(124,20){\vector(0,1){1.00}}

\put(110,10){\vector(-1,0){1.00}}
\put(110,30){\vector(1,0){1.00}}
\put(100,20){\vector(0,1){1.00}}
\put(120,20){\vector(0,-1){1.00}}

\put(108.0,11.3){$\be_4$}
\put(108.0,27.0){$\be_2$}
\put(101,19.0){$\be_1$}
\put(115.5,19.0){$\be_3$}

\put(108.0,3.0){$\be_8$}
\put(108.0,35.5){$\be_6$}
\put(91.5,19.0){$\be_5$}
\put(125,19.0){$\be_7$}
\end{picture}

\caption{\emph{a}) The square lattice $\dL^2$;
the fundamental cell $\Omega$ is shaded; \;  \emph{b}) the fundamental graph $\dL_*^2$.} \label{slex}
\end{figure}
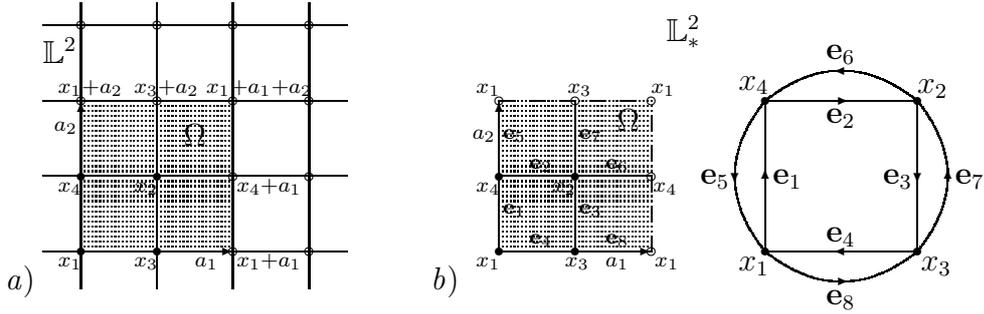

\begin{example}\lb{ExSL}
Let $H=A-\vk+V$ be the Schr\"odinger operator with a periodic
potential $V$  on the square lattice $\dL^2$. We assume that the
potential $V$ has the periods $a_1$, $a_2$, see Fig.~\ref{slex}a.
Then the traces of $A^n(k)$ and $H^n(k)$, $k=(k_1,k_2)\in\T^2$,
given by \er{TrAo} and \er{TrnH}, for $n=1,2,3$ have the form
\[\lb{TASL}
\begin{array}{ll}
\Tr A(k)=0, & \Tr H(k)=-16+\sum\limits_{s=1}^4V_s, \\
\Tr A^2(k)=8\cF(k), \qqq & \Tr H^2(k)=\sum\limits_{s=1}^4v_s^2+\Tr A^2(k),  \\
\Tr A^3(k)=0, & \Tr H^3(k)=\sum\limits_{s=1}^4\big(v_s^3+6\cF(k)v_s\big),
\end{array}
\]
where
\[\lb{Ex1F}
\cF(k)=2+\cos k_1+\cos k_2, \qqq v_s=V_s-4, \qqq V_s=V_{x_s},\qqq s=1,2,3,4.
\]
\end{example}

\section{Proof of the main results}
\setcounter{equation}{0}
\lb{Sec3}

\subsection{Trace formulas for adjacency operators}
First we prove properties of the numbers $\cN_n$ and $\cN_n^\mm$.

\medskip

\no \textbf{Proof of Proposition \ref{spcN}.} \emph{i}) Recall that
$\cN_n=\#\cC_n$ and $\cN_n^\mm=\#\cC_n^\mm$. Since there are no
loops in the periodic graph $\cG$, the fundamental graph $\cG_*$ has
no loops with zero index. This yields the first identity in
\er{cN10}.

Let $n\in\N$. Each oriented edge $\be\in\cA_*$ generates the cycle $\bc_\be=(\be,\ul\be,\ldots,\be,\ul\be)$ of length $2n$ consisting of $n$ backtracks. The index of the cycle $\bc_\be$ is 0. Thus, the number $\cN_{2n}^0$ of all cycles of length $2n$ with zero index is not less than the number of edges from $\cA_*$. If there are no multiple edges in the periodic
graph, then all cycles of length $2$ with zero index in $\cG_*$ are
cycles $\bc_\be=(\be,\ul\be)$, $\be\in\cA_*$, and $\cN_2^0=\#\cA_*$.

We show the first inequality in \er{Npme}. We fix a vertex $x\in\cV_*$ and estimate the number of cycles $\bc$ of length $n$ starting at $x$. Since each edge of $\bc$ may be chosen in not more than $\vk_+$ ways, then the cycle $\bc$ may be chosen in not more than $\vk_+^n$ ways. Repeating this for each vertex $x\in\cV_*$, we obtain the first inequality in \er{Npme}.

For each cycle $\bc=(\be_1,\ldots,\be_n)\in\cC_n$ there exists a
cycle $\ul{\bc}=(\ul{\be}_n,\ldots,\ul{\be}_1)\in\cC_n$ and, by
\er{ininc},  $\t(\ul{\bc}\,)=-\t(\bc)$. This yields the second identity in \er{Npme}.

By the definition of the cycle index \er{cyin}, for each cycle $\bc$
of length $n$ we have
\[\lb{tpl}
\textstyle\|\t(\bc)\|\leq \sum\limits_{\be\in\bc}\|\t(\be)\|\leq
n\t_+.
\]
Thus, we obtain \er{cNe0}.

\emph{ii}) Let $\cG_*$ be bipartite. Then it has no odd-length cycles. Thus, the identities \er{cNBG} hold true. Let $\cN_n=0$ for all odd $n\leq\n$. Then there are no cycles of odd length $\leq\n$ in $\cG_*$. Assume that there exists a cycle $\bc$ of odd length $>\n$ in $\cG_*$. Since $\cG_*$ is a graph with $\n$ vertices, the cycle $\bc$ has repeated vertices (other than the first and last ones). Then there exists a cycle of odd length $\leq\n$ in $\cG_*$. We get a contradiction. Thus, there are no cycles of any odd length in $\cG_*$. This yields that $\cG_*$ is bipartite.

\emph{iii}) Let $\cG$ be bipartite. Then there are no odd-length cycles in $\cG$. Consequently, the fundamental graph $\cG_*$ has no odd-length cycles with zero index. (Recall that each cycle with zero index in the fundamental graph $\cG_*$ corresponds to a cycle of the same length in the periodic graph $\cG$, see Remark 1 in Subsection \ref{SoC}.) Thus, $\cN_n^0=0$ for odd $n\in\N$. Conversely, let $\cN_n^0=0$ for all odd $n$, i.e., there are no odd-length cycles with zero index in $\cG_*$. This yields that there are no odd-length cycles in $\cG$, and therefore $\cG$ is bipartite.

\emph{iv}) The last statement is proved in Example \ref{ExBB}.
\qq \BBox

\medskip

From Theorem \ref{TFDC} it follows that in the standard orthonormal
basis of $\ell^2(\cV_*)=\C^\n$, $\n=\#\cV_*$, the $\n\ts\n$ matrix
$A(k)=\big(A_{xy}(k)\big)_{x,y\in\cV_*}$ of the fiber adjacency
operator $A(k)$ given by \er{fado} has the form
\[\lb{fnml0}
A_{xy}(k)=\sum_{\be=(x,y)\in\cA_*}e^{-i\lan\t(\be),k\ran}.
\]
The diagonal entries of the matrix $A(k)$ may also be written in the form
\[\lb{deDa}
A_{xx}(k)=\sum\limits_{\be=(x,x)\in\cA_*}
\cos\lan\t(\be),k\ran, \qqq \forall\,x\in\cV_*.
\]

Now we prove Theorem \ref{TFao} about trace formulas for the fiber
adjacency operator $A(k)$. We need the following simple identity
\[\lb{inex}
\frac1{(2\pi)^d}\int_{\T^d}\cos\lan\mm,k\ran\,dk=\d_{\mm,0}, \qqq \forall\,\mm\in\Z^d,
\]
where $\d_{\mm,0}$ is the Kronecker delta.

\medskip

\no\textbf{Proof of Theorem \ref{TFao}.} \emph{i})   Using
\er{fnml0}, for each $n\in\N$ we obtain
\begin{multline*}
\Tr A^n(k)=\sum_{x_1,\ldots,x_n\in\cV_*}A_{x_1x_2}(k)A_{x_2x_3}(k)
\ldots A_{x_{n-1}x_n}(k)A_{x_nx_1}(k)
\\
=
\sum_{x_1,\ldots,x_n\in\cV_*}\sum\limits_{\be_1,\ldots,\be_n\in\cA_*}
e^{-i\lan\t(\be_1)+\t(\be_2)+\ldots+\t(\be_n),k\ran}
=\sum_{\bc\in\cC_n}e^{-i\lan\t(\bc),k\ran},
\end{multline*}
where $\be_j=(x_j,x_{j+1})$, $j\in \N_n$, $x_{n+1}=x_1$. Here we
have also used the definition \er{cyin} of the cycle indices. Thus, we
have the finite Fourier series for the $2\pi\Z^d$-periodic function $\Tr
A^n(k)$, since $\t(\bc)\in \Z^d$. We rewrite this Fourier series in
the standard form
\[\label{rg11+}
\Tr A^n(k)=\sum_{\bc\in\cC_n}e^{-i\lan\t(\bc),k\ran}
=\sum\limits_{\mm\in\Z^d}\sum_{\bc\in\cC_n^\mm}e^{-i\lan\mm,k\ran}
=\sum\limits_{\mm\in\Z^d}e^{-i\lan\mm,k\ran}\cN_n^\mm,\qqq \cN_n^\mm=\#\cC_n^\mm,
\]
or, using \er{Npme} and \er{cNe0},
\[\lb{Trak}
\Tr A^n(k)=\sum_{\bc\in\cC_n}\cos\lan\t(\bc),k\ran
=\sum\limits_{\mm\in\Z^d}\cN_n^\mm\cos\lan\mm,k\ran=
\sum_{\mm\in\Z^d\atop\|\mm\|\leq n\t_+}\cN_n^\mm\cos\lan\mm,k\ran.
\]
Using \er{Trak}, we obtain
$$
\big|\Tr A^n(k)\big|=\bigg|\sum_{\bc\in\cC_n}\cos\lan\t(\bc),k\ran\bigg|\leq
\Tr A^n(0)=\sum\limits_{\mm\in\Z^d}\cN_n^\mm=\cN_n.
$$

\emph{ii}) Integrating \er{TrAo} over $k\in\T^d$ and using
\er{inex},  we obtain \er{ITrA}.

\emph{iii}) Let $\cG_*$ be bipartite. Then, by \er{cNBG},
$\cN_n^\mm=0$  for odd $n\in\N$ and all $\mm\in\Z^d$. This and
\er{TrAo} yield \er{Kbfg}. Conversely, let the condition \er{Kbfg}
be fulfilled. Then, by \er{EsTAn}, $\cN_n=0$ for all odd $n\leq\n$, and, due to Proposition \ref{spcN}.\emph{ii}, $\cG_*$ is bipartite.

\emph{iv}) By Proposition \ref{spcN}.\emph{iii}, $\cG$ is bipartite iff $\cN_n^0=0$ for odd $n\in\N$. This and \er{ITrA} prove the first statement.
Due to Proposition \ref{spcN}.\emph{iv}, for any $s\in\N$ there exists a periodic graph $\cG$ such that
$$
\cN_n^0=0\qq \textrm{for all odd $n<2s$,\qq and} \qq
\cN_{2s+1}^0\neq0.
$$
This and \er{ITrA} yield that $\int_{\T^d}\Tr A^n(k)dk=0$ for all odd $n<2s$ and there exists a cycle of odd length $2s+1$ with zero index in $\cG_*$. This cycle corresponds to a cycle in $\cG$ with the same odd length $2s+1$. Thus, $\cG$ is non-bipartite. \qq $\BBox$

\begin{remark} Theorem \ref{TFao}.\emph{i} gives another simple proof of the first inequality in \er{Npme}. Since $\s(A)\subseteq[-\vk_+,\vk_+]$, then $\s\big(A(k)\big)\ss[-\vk_+,\vk_+]$ for all $k\in\T^d$ and, using \er{TrAo} and \er{EsTAn}, we have
$$
\cN_n=\Tr A^n(0)=\sum_{j=1}^\n\big(\l_j^o(0)\big)^n\leq
\sum_{j=1}^\n\big|\l_j^o(0)\big|^n\leq\n\vk_+^n.
$$
\end{remark}

The following example shows that the condition $\cN_n^0=0$ for all
odd $n\leq\n$ is not enough for a periodic graph to be bipartite.

\begin{example}\lb{ExBB}
Let $p\in\N$. Consider a $\Z$-periodic graph $\cG_p$, obtained from
the  one-dimensional lattice $\Z$ by adding the edges $(n,n+p)$, $n\in\N$
(Fig.\ref{figB}a). Let $A$ be the adjacency operator on $\cG_p$.
Then the following statements hold true.

$\bu$ The fundamental graph $\cG_*=\cG_p/\Z$ is non-bipartite for any $p\in\N$.

$\bu$ If $p$ is odd, then the condition \er{cN0BG} is fulfilled, and $\cG_p$ is bipartite.

$\bu$ If $p$ is even, then $\cN_n^0=0$ for odd $n=1,3,\ldots,p-1$, but
$\cN_{p+1}^0\neq0$ and $\cG_p$ is non-bipartite.
\end{example}

\setlength{\unitlength}{1.0mm}
\begin{figure}[h]
\centering
\unitlength 1.0mm 
\begin{picture}(160,25)
\put(5,10){\line(1,0){100.00}} \put(10,10){\circle*{1}}
\put(25,10){\circle*{1}} \put(40,10){\circle*{1}}
\put(55,10){\circle*{1}} \put(70,10){\circle*{1}}
\put(85,10){\circle*{1}} \put(100,10){\circle*{1}}
\put(18,10){\vector(1,0){1}} \put(40,20){\vector(1,0){1}}
\put(16,7){$\be_1$} \put(39,21.5){$\be_2$}
\bezier{600}(10,10)(40,30)(70,10) \bezier{600}(25,10)(55,30)(85,10)
\bezier{600}(40,10)(70,30)(100,10)
\bezier{600}(55,10)(85,28)(105,14)
\bezier{600}(70,10)(90,22)(105,20)
\bezier{600}(85,10)(98,18)(105,19)
\bezier{600}(100,10)(102,12)(105,13)

\bezier{600}(55,10)(25,28)(5,14) \bezier{600}(40,10)(20,22)(5,20)
\bezier{600}(25,10)(12,18)(5,19) \bezier{600}(10,10)(8,12)(5,13)
\put(0,6){\emph{a})} \put(0,22){$\cG_p$} \put(9,6){$0$}
\put(24,6){$1$} \put(39,6){$2$} \put(54,6){$3$} \put(60,6){$\ldots$}
\put(69,6){$p$} \put(81,6){$p+1$} \put(96,6){$p+2$}
\put(114,6){\emph{b})} \put(113,22){$\cG_*$} \put(116.5,14){$0$}
\put(130,15){\vector(0,-1){1}} \put(140,15){\vector(0,-1){1}}
\put(131,14){$\be_1$} \put(141,14){$\be_2$}
\put(120,15){\circle*{1}} \bezier{600}(120,15)(121,24)(130,25)
\bezier{600}(140,15)(139,24)(130,25)
\bezier{600}(120,15)(121,6)(130,5)
\bezier{600}(140,15)(139,6)(130,5)

\bezier{600}(120,15)(120.5,19.5)(125,20)
\bezier{600}(130,15)(129.5,19.5)(125,20)
\bezier{600}(120,15)(120.5,10.5)(125,10)
\bezier{600}(130,15)(129.5,10.5)(125,10)

\end{picture}
\vspace{-4mm} \caption{\footnotesize \emph{a}) A $\Z$-periodic graph
$\cG_p$, $p\in\N$; if $p$ is odd, then $\cG_p$ is bipartite; if $p$
is even, then $\cG_p$ is non-bipartite; \;  \emph{b}) the
fundamental graph $\cG_*=\cG_p/\Z$ consists of one vertex 0 and two
loop edges $\be_1,\be_2$ with indices $\t(\be_1)=1$, $\t(\be_2)=p$;
$\cG_*$ is non-bipartite for any $p\in\N$.} \label{figB}
\end{figure}
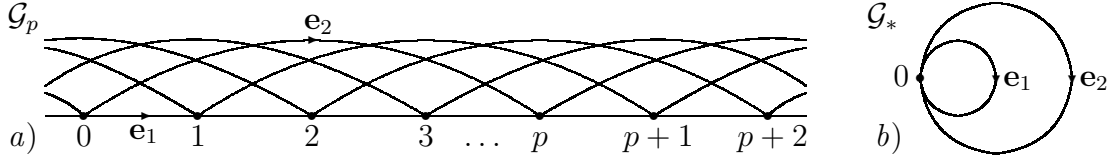

\no \textbf{Proof.} The fundamental graph $\cG_*=\cG_p/\Z$ of $\cG_p$
consists of one vertex 0 and two loop edges $\be_1,\be_2$ with
indices $\t(\be_1)=1$, $\t(\be_2)=p$, see Fig.\ref{figB}\emph{b}.
Since $\cG_*$ has loops, it is non-bipartite for any $p\in\N$. If
$p$ is odd, then there are no odd-length cycles with zero index in
$\cG_*$, i.e., the condition \er{cN0BG} is fulfilled, and, consequently, the periodic graph $\cG_p$ is
bipartite. If $p$ is even, then there are still no cycles with zero
index of length $1,3,\ldots,p-1$ in $\cG_*$. But the cycle
$(\underbrace{\be_1,\ldots\be_1}_{p\textrm{ times}},\ul\be_2\,)$ has
zero index and the odd length $p+1$. Thus, $\cN_n^0=0$ for odd $n=1,3,\ldots,p-1$, and $\cN_{p+1}^0\neq0$. In the case of even $p$ the periodic graph $\cG_p$ is non-bipartite, since there is a cycle of odd length $p+1$ in $\cG_p$. \qq \BBox

\subsection{Trace formulas for Schr\"odinger operators}
First we prove Proposition \ref{Ppof} about properties of the
functionals  $\cT_n(k)$ defined by \er{deTnk}.

\medskip

\no\textbf{Proof of Proposition \ref{Ppof}.} \emph{i}) Using
\er{Wcy}  and \er{obc+}, we have $
\big|\o(\bc)\big|\leq\o_+^{|\bc|}$ for all $\bc\in\wt\cC$, which
yields
$$
\big|\cT_n(k)\big|=\bigg|\sum_{\bc\in\wt\cC_n}\o(\bc)\cos\lan\t(\bc),k\ran\bigg|\leq
\sum_{\bc\in\wt\cC_n}\big|\o(\bc)\big|=\sum_{\bc\in\wt\cC_n}\o_+^n=\o_+^n\,\wt\cN_n.
$$

\emph{ii}) The identity \er{deTnk} is a finite Fourier series for
the $2\pi\Z^d$-periodic function $\cT_n(k)$, since $\t(\bc)\in\Z^d$. We
rewrite this Fourier series in the standard form
\[\lb{rg11H}
\begin{aligned}
\cT_n(k)=\sum_{\bc\in\wt\cC_n}\o(\bc)\cos\lan\t(\bc),k\ran
=\sum_{\mm\in\Z^d}\sum_{\bc\in\wt\cC_n^\mm}
\o(\bc)\cos\lan\mm,k\ran
\\
=\sum\limits_{\mm\in\Z^d}\cos\lan\mm,k\ran\sum_{\bc\in\wt\cC_n^\mm}\o(\bc)
=\sum\limits_{\mm\in\Z^d}\cT_{n,\mm}\cos\lan\mm,k\ran,
\end{aligned}
\]
where the coefficients $\cT_{n,\mm}$ are given in \er{FsTrH}. Due to
\er{cNe0}, $\cT_{n,\mm}=0$ for all $\mm\in\Z^d$ such that
$\|\mm\|>n\t_+$, i.e., the sum in RHS of \er{rg11H} is finite. Thus,
the identity \er{rg11H} can be written in the form \er{FsTrH}.
Integrating \er{FsTrH} over $k\in\T^d$ and using \er{inex}, we
obtain \er{IVFu}. \qq \BBox

\medskip

From Theorem \ref{TFDC} it follows that in the standard orthonormal
basis of  $\ell^2(\cV_*)=\C^\n$, $\n=\#\cV_*$, the $\n\ts\n$ matrix
$H(k)=\big(H_{xy}(k)\big)_{x,y\in\cV_*}$ of the fiber Schr\"odinger
operator $H(k)$ defined by \er{Hvt'} -- \er{fado} has the form
\[\lb{mrSo}
H(k)=A(k)+v, \qqq\textrm{where}\qqq v=\diag(v_x)_{x\in\cV_*}, \qqq v_x=V_x-\vk_x,
\]
and the $\n\ts\n$ matrix $A(k)$ of the fiber adjacency operator is given by \er{fnml0}.

The next lemma shows that the operator $H(k)$ can be considered as a
fiber weighted adjacency operator on the modified fundamental graph $\wt\cG_*=(\cV_*,\wt\cA_*)$.

\begin{lemma}\lb{LDeH}
The fiber Schr\"odinger operator $H(k)=\big(H_{xy}(k)\big)_{x,y\in\cV_*}$
given by \er{mrSo}, \er{fnml0} satisfies
\[\lb{mrmL+}
H_{xy}(k)=\sum_{\be=(x,y)\in\wt\cA_*}\o(\be) e^{-i\lan\t(\be),\,k\ran},\qqq
\forall\, x,y\in\cV_*, \qqq \forall\,k\in\T^d,
\]
where $\wt\cA_*$ is the set of all edges of the modified fundamental
graph $\wt\cG_*$ defined by \er{wtAs}; $\o(\be)$ is given by
\er{webe}, and $\t(\be)$  is the index of the edge $\be\in\wt\cA_*$
defined by \er{in}, \er{dco}.
\end{lemma}

\no \textbf{Proof.} Let $x,y\in\cV_*$. If $x\neq y$, then, using
\er{webe}  and \er{mrSo}, \er{fnml0}, we have
$$
\sum_{\be=(x,y)\in\wt\cA_*}\o(\be) e^{-i\lan\t(\be),\,k\ran}
=\sum_{\be=(x,y)\in\cA_*} e^{-i\lan\t(\be),\,k\ran}=H_{xy}(k).
$$
Similarly, if $x=y$, then we obtain
$$
\sum_{\be=(x,x)\in\wt\cA_*}\o(\be) e^{-i\lan\t(\be),\,k\ran}
=\o(\be_x)+\sum_{\be=(x,x)\in\cA_*} e^{-i\lan\t(\be),\,k\ran}
=v_x+\sum_{\be=(x,x)\in\cA_*} e^{-i\lan\t(\be),\,k\ran}=H_{xx}(k).
$$
Thus, the identity \er{mrmL+} has been proved. \qq \BBox

\medskip

Now we prove Theorems \ref{TPG} and \ref{TPG2} about trace formulas
for  the fiber Schr\"odinger operator $H(k)$.

\medskip

\no\textbf{Proof of Theorem \ref{TPG}.} The proof is similar  to the
proof  of Theorem \ref{TFao}.\emph{i}. Using \er{mrmL+} and \er{deTnk}, for
each $n\in\N$ we obtain
\begin{multline*}
\Tr H^n(k)=
\sum_{x_1,\ldots,x_n\in\cV_*}H_{x_1x_2}(k)H_{x_2x_3}(k)\ldots H_{x_{n-1}x_n}(k)H_{x_nx_1}(k)
\\
=\!\!\!\! \sum_{x_1,\ldots,x_n\in\cV_*}
\sum_{\be_1,\ldots,\be_n\in\wt\cA_*}\!\!\o(\be_1)\o(\be_2)\ldots\o(\be_n)
e^{-i\lan\t(\be_1)+\t(\be_2)+\ldots+\t(\be_n),k\ran} =
\sum_{\bc\in\wt\cC_n}\o(\bc)e^{-i\lan\t(\bc),k\ran}=\cT_n(k),
\end{multline*}
where $\be_j=(x_j,x_{j+1})$, $j\in \N_n$, and $x_{n+1}=x_1$. Here
we have also used the definitions \er{cyin} and \er{Wcy} of the cycle
index $\t(\bc)$ and the cycle weight $\o(\bc)$. Thus, \er{TrnH} is
proved.

Now we prove \er{TrH123}. Using \er{deTnk} and \er{oc1C},
the identity \er{TrnH} for $n=1$ has the form
\[\lb{THco1}
\Tr H(k)=\cT_1(k)=\sum_{\bc\in\wt\cC_1}\o(\bc)\cos\lan\t(\bc),k\ran=
\sum_{\bc\in\wt\cC_1\sm\cC_1}\o(\bc)\cos\lan\t(\bc),k\ran+
\sum_{\bc\in\cC_1}\cos\lan\t(\bc),k\ran.
\]
Since $\bc\in\wt\cC_1\sm\cC_1$ iff $\bc$ is an added loop $\be_x$ at
some vertex $x\in\cV_*$ with index $\t(\be_x)=0$ and weight
$\o(\be_x)=v_x$, we obtain the first identity in \er{TrH123}. In
particular, if there are no loops in the fundamental graph $\cG_*$,
i.e., $\cC_1=\varnothing$, then this identity has the form \er{hnl}.
Similarly, the identity \er{TrnH} for $n=2$  has the form
\[\lb{THco2}
\Tr H^2(k)=\cT_2(k)=\sum_{\bc\in\wt\cC_2}\o(\bc)\cos\lan\t(\bc),k\ran=
\sum_{\bc\in\wt\cC_2\sm\cC_2}\o(\bc)\cos\lan\t(\bc),k\ran+
\sum_{\bc\in\cC_2}\cos\lan\t(\bc),k\ran.
\]
Each cycle from $\wt\cC_2\sm\cC_2$ has one of the following forms:

$\bu$ $\bc_1=(\be_x,\be_x)$ for some $x\in\cV_*$, where $\be_x$ is
a loop at the vertex $x$ with index $\t(\be_x)=0$ and weight
$\o(\be_x)=v_x$ added to the fundamental graph $\cG_*$, see
\er{wtAs};

$\bu$ $\bc_2=(\be_x,\bc)$ or $\bc_3=(\bc,\be_x)$, where $\bc\in\cC_1$
is a loop at some vertex $x$ (i.e., a cycle of length one) in $\cG_*$ with weight $\o(\bc)=1$.
\\
Due to the definitions \er{cyin} and \er{Wcy} of the cycle index and the cycle weight, we have
$$
\begin{array}{ll}
\t(\bc_1)=0,\qqq & \t(\bc_2)=\t(\bc_3)=\t(\bc),\\[4pt]
\o(\bc_1)=v_x^2,\qq & \o(\bc_2)=\o(\bc_3)=v_x=v(\bc).
\end{array}
$$
Then the first sum in RHS of \er{THco2} has the form
\[\lb{fsu}
\sum_{\bc\in\wt\cC_2\sm\cC_2}\o(\bc)\cos\lan\t(\bc),k\ran=
\sum_{x\in\cV_*}v_x^2+2\sum_{\bc\in\cC_1}v(\bc)\cos\lan\t(\bc),k\ran.
\]
Combining \er{THco2} and \er{fsu}, we obtain the second identity in
\er{TrH123}.  If there are no loops and multiple edges in the
fundamental graph $\cG_*$, then $\cC_1=\varnothing$ and all cycles
of length 2 in $\cG_*$ have the form $\bc=(\be,\ul\be\,)$, where
$\be\in\cA_*$, and $\t(\bc)=0$. This yields
$$
\sum\limits_{\bc\in\cC_2}\cos\lan\t(\bc),k\ran=\#\cA_*,
$$
and, using the well-known identity
$\sum\limits_{x\in\cV_*}\vk_x=\#\cA_*$,  we obtain \er{nome}. \qq
$\BBox$

\medskip

\no\textbf{Proof of Theorem \ref{TPG2}.} Integrating \er{TrnH}  over
$k\in\T^d$ and using \er{IVFu}, we obtain \er{ITrH}. Similarly,
integrating \er{TrH123} over $k\in\T^d$ and using \er{inex} and the identity $\cN_1^0=0$, we get
\er{TR13IH}. If there are no multiple edges in the periodic graph
$\cG$, then, by \er{cN0g1}, $\cN_2^0=\#\cA_*$. \qq \BBox

\medskip

In the next proposition we present the third order traces for the  Schr\"odinger operator.

\begin{proposition}\lb{TrTO}
Let $H(k)=-\D(k)+V$, $k\in\T^d$, be the fiber Schr\"odinger operator defined
by \er{Hvt'} -- \er{fado} on the fundamental graph $\cG_*=(\cV_*,\cA_*)$. Then
\begin{multline}\label{todT}
\Tr H^3(k)=\cT_3(k)=\sum_{\bc\in\cC_3}\cos\lan\t(\bc),k\ran+
\sum_{x\in\cV_*}v_x^3\\+3\sum_{\bc\in\cC_1}v^2(\bc)\cos\lan\t(\bc),k\ran
+\frac32\sum_{\bc\in\cC_2}v(\bc)\cos\lan\t(\bc),k\ran, \qqq v_x=V_x-\vk_x,
\end{multline}
\[\lb{todTI}
\frac1{(2\pi)^d}\int_{\T^d}\Tr H^3(k)dk=\cT_{3,0}=\cN_3^0+
\sum_{x\in\cV_*}v_x^3+
\frac32\sum_{\bc\in\cC_2^0}v(\bc).
\]
If there are no multiple edges in the periodic graph $\cG$, then
\[\lb{HA3n}
\frac1{(2\pi)^d}\int_{\T^d}\Tr H^3(k)dk=\cN_3^0+\sum_{x\in\cV_*}\big(v_x^3+3v_x\vk_x\big).
\]
\end{proposition}

\no \textbf{Proof.} Using \er{deTnk} and \er{oc1C}, the identity \er{TrnH} for $n=3$ has the form
\[\lb{THco3}
\Tr H^3(k)=\cT_3(k)=\sum_{\bc\in\wt\cC_3}\o(\bc)\cos\lan\t(\bc),k\ran=
\sum_{\bc\in\wt\cC_3\sm\cC_3}\o(\bc)\cos\lan\t(\bc),k\ran+
\sum_{\bc\in\cC_3}\cos\lan\t(\bc),k\ran.
\]
Each cycle from $\wt\cC_3\sm\cC_3$ has one of the following forms

$\bu$ $\bc_1=(\be_x,\be_x,\be_x)$ for some $x\in\cV_*$, where
$\be_x$ is  a loop at the vertex $x$ with index $\t(\be_x)=0$ and
weight $\o(\be_x)=v_x$ added to the fundamental graph $\cG_*$, see
\er{wtAs};

$\bu$ $\bc_2=(\be_x,\be_x,\bc)$ or $\bc_2=(\be_x,\bc,\be_x)$ or
$\bc_2=(\bc,\be_x,\be_x)$, where $\bc\in\cC_1$ is a loop at some
vertex  $x$ (i.e., a cycle of length one) in $\cG_*$;

$\bu$ $\bc_3=(\be_x,\be_1,\be_2)$ or $\bc_3=(\be_1,\be_2,\be_x)$ or
$\bc_4=(\be_1,\be_y,\be_2)$, where $\be_1=(x,y)$, $\be_2=(y,x)$
for some vertices $x,y\in\cV_*$, i.e., $\wt\bc:=(\be_1,\be_2)\in\cC_2$ is a cycle of length two in $\cG_*$.\\
Due to the definitions \er{cyin} and \er{Wcy} of the cycle index and the cycle weight, we have
$$
\begin{array}{llll}
\t(\bc_1)=0,\qqq & \t(\bc_2)=\t(\bc),\qqq & \t(\bc_3)=\t(\wt\bc\,),\qqq & \t(\bc_4)=\t(\wt\bc\,),\\[4pt]
\o(\bc_1)=v_x^3,\qqq & \o(\bc_2)=v_x^2=v^2(\bc), \qqq &
\o(\bc_3)=v_x, \qqq &
\o(\bc_4)=v_y.
\end{array}
$$
Then the first sum in RHS of \er{THco3} has the form
\[\label{fsu3}
\begin{aligned}
&\sum_{\bc\in\wt\cC_3\sm\cC_3}\o(\bc)\cos\lan\t(\bc),k\ran=
\sum_{x\in\cV_*}v_x^3+3\sum_{\bc\in\cC_1}v^2(\bc)\cos\lan\t(\bc),k\ran+
\sum_{\wt\bc\in\cC_2}(2v_x+v_y)\cos\lan\t(\wt\bc\,),k\ran,\\
&\hspace{20mm}\textrm{where}\qqq \wt\bc=(\be_1,\be_2),\qqq
\be_1=(x,y),\qqq \be_2=(y,x).
\end{aligned}
\]
We rewrite the last sum in the form
\begin{multline}\lb{fsu33}
\sum_{\wt\bc=(\be_1,\be_2)\in\cC_2}(2v_x+v_y)\cos\lan\t(\wt\bc\,),k\ran\\=
\frac12\sum_{\wt\bc=(\be_1,\be_2)\in\cC_2}(2v_x+v_y)\cos\lan\t(\wt\bc\,),k\ran+
\frac12\sum_{\wt\bc=(\be_2,\be_1)\in\cC_2}(2v_y+v_x)\cos\lan\t(\wt\bc\,),k\ran\\
=\frac32\sum_{\wt\bc=(\be_1,\be_2)\in\cC_2}(v_x+v_y)\cos\lan\t(\wt\bc\,),k\ran=
\frac32\sum_{\wt\bc\in\cC_2}v(\wt\bc\,)\cos\lan\t(\wt\bc\,),k\ran.
\end{multline}
Combining \er{THco3} -- \er{fsu33}, we obtain \er{todT}. Integrating \er{todT}
 over $k\in\T^d$ and using \er{inex} and \er{cN10}, we get \er{todTI}.

If there are no multiple edges in the periodic graph $\cG$, then all
cycles  of length 2 and with zero index in $\cG_*$ have the form
$\bc=(\be,\ul\be\,)$, where $\be\in\cA_*$. This yields
$$
\sum_{\bc\in\cC_2^0}v(\bc)=\sum_{\be=(x,y)\in\cA_*}(v_x+v_y)=
2\sum_{x\in\cV_*}\vk_xv_x,
$$
and the identity \er{todTI} has the form \er{HA3n} \qq \BBox

\medskip

Now we prove Corollary \ref{CCrFA} about the difference of  the
traces of the fiber Schr\"odinger operators and the corresponding
fiber adjacency operators.

\medskip

\no\textbf{Proof of Corollary \ref{CCrFA}.}  Due to \er{TrAo},
\er{TrnH}  and \er{deTnk}, \er{FsTrH}, we have
\begin{multline}\label{difT}
\Tr\big(H^n(k)-A^n(k)\big)=\sum_{\bc\in\wt\cC_n}\o(\bc)\cos\lan\t(\bc),k\ran-
\sum_{\bc\in\cC_n}\cos\lan\t(\bc),k\ran\\=\sum\limits_{\mm\in\Z^d\atop
\|\mm\|\leq n\t_+}(\cT_{n,\mm}-\cN_n^\mm)\cos\lan\mm,k\ran.
\end{multline}
Using the inclusion $\cC_n\ss\wt\cC_n$, the identity \er{oc1C} and
the definition of $\cT_{n,\mm}$ in \er{FsTrH}, we can rewrite
\er{difT} in the form \er{TrHAn}.

Integrating \er{TrHAn} over $k\in\T^d$ and using \er{inex}, we get
\er{TrHAnI}. Combining \er{TrAo} (as $n=1,2$) and \er{TrH123}, we
obtain \er{TrHA1}.  Similarly, combining \er{ITrA} and \er{TR13IH}
(or integrating \er{TrHA1} over $k\in\T^d$), we have \er{TrHA1I}.\qq
\BBox

\medskip

In the following statement we compare the traces of the fiber
Schr\"odinger  operators and the traces of the corresponding fiber
Laplacians.

\begin{corollary}\lb{CHmL}
Let $H(k)=H_0(k)+V$, $H_0=-\D(k)$, $k\in\T^d$, be the fiber
Schr\"odinger  operator defined by \er{Hvt'} -- \er{fado}. Then for
each $n\in\N$
\begin{multline}\label{TrHDn}
\Tr\big(H^n(k)-H_0^n(k)\big)
=\sum_{\bc\in\wt\cC_n\sm\cC_n}\big(\o(\bc,V)-\o(\bc,0)\big)\cos\lan\t(\bc),k\ran\\
=\sum_{\mm\in\Z^d\atop\|\mm\|\leq n\t_+}\gt^o_{n,\mm}\cos\lan\mm,k\ran, \qqq \gt^o_{n,\mm}=
\sum_{\bc\in\wt\cC_n^\mm\sm\cC_n^\mm}\big(\o(\bc,V)-\o(\bc,0)\big),
\end{multline}
\[\label{TrHDnI}
\frac1{(2\pi)^d}\int_{\T^d}\Tr\big(H^n(k)-H_0^n(k)\big)dk
=\gt^o_{n,0}, \qqq \gt^o_{n,0}=
\sum_{\bc\in\wt\cC_n^0\sm\cC_n^0}\big(\o(\bc,V)-\o(\bc,0)\big),
\]
where $\t_+$ is defined in \er{cNe0}, and $\|\cdot\|$ is the
standard norm in $\R^d$.  In particular,
\[\label{TrHD1}
\begin{aligned}
&\Tr\big(H(k)-H_0(k)\big)=\sum_{x\in\cV_*}V_x,\\
&\Tr\big(H^2(k)-H_0^2(k)\big)=\sum_{x\in\cV_*}\big(V_x^2-2V_x\vk_x\big)+
2\sum_{\bc\in\cC_1}V(\bc)\cos\lan\t(\bc),k\ran,\\
&\Tr\big(H^3(k)-H_0^3(k)\big)=\sum_{x\in\cV_*}\big(V_x^3-3V_x^2\vk_x+3V_x\vk_x^2\big)
\\&\hspace{36mm}+3\sum_{\bc\in\cC_1}\big(V^2(\bc)-2V(\bc)\vk(\bc)\big)\cos\lan\t(\bc),k\ran
+\frac32\sum_{\bc\in\cC_2}V(\bc)\cos\lan\t(\bc),k\ran,
\end{aligned}
\]
\[\label{TrHD1I}
\begin{aligned}
&\frac{1}{(2\pi)^d}\int_{\T^d}\Tr\big(H(k)-H_0(k)\big)dk=\sum_{x\in\cV_*}V_x,\\
&\frac1{(2\pi)^d}\int_{\T^d}\Tr\big(H^2(k)-H_0^2(k)\big)dk=
\sum_{x\in\cV_*}\big(V_x^2-2V_x\vk_x\big),\\
&\frac{1}{(2\pi)^d}\int_{\T^d}\Tr\big(H^3(k)-H_0^3(k)\big)dk=
\sum_{x\in\cV_*}\big(V_x^3-3V_x^2\vk_x+3V_x\vk_x^2\big)+
\frac32\sum_{\bc\in\cC_2^0}V(\bc),
\end{aligned}
\]
where
$$
\begin{aligned}
&V(\bc)=V_{x_1}+\ldots+V_{x_n},\qqq \vk(\bc)=\vk_{v_1}+\ldots+\vk_{v_n},
\\[2pt]
& \textrm{for each}\qq \bc=(\be_1,\ldots,\be_n)\in\cC_n,
\qq  \be_j=(x_j,x_{j+1}), \qq j\in \N_n, \qq x_{n+1}=x_1.
\end{aligned}
$$
\end{corollary}

\no \textbf{Proof.} Using the identity
$$
\Tr\big(H^n(k)-H^n_0(k)\big)=\Tr\big(H^n(k)-A^n(k)\big)-\Tr\big(H_0^n(k)-A^n(k)\big),
$$
and applying to its RHS Corollary \ref{CCrFA} and Proposition
\ref{TrTO},  we obtain the required statements. \qq \BBox

\subsection{Trace formulas for the heat kernel and the resolvent}
Now we prove Corollaries \ref{CHC} and \ref{CHC1} about trace
formulas for  the heat kernel and the resolvent of the fiber
Schr\"odinger operators.

\medskip

\no\textbf{Proof of Corollary \ref{CHC}.} Using \er{TrnH} and \er{deTnk}, we have
\begin{multline*}
\Tr e^{t H(k)}=\Tr\sum_{n=0}^\iy\frac{t^n}{n!}\,H^n(k)=
\sum_{n=0}^\iy\frac{t^n}{n!}\,\Tr H^n(k)=\n+\sum_{n=1}^\iy\frac{t^n}{n!}\,\cT_n(k)\\=
\n+\sum_{n=1}^\iy\frac{t^n}{n!}\sum_{\bc\in\wt\cC_n}\o(\bc)\cos\lan\t(\bc),k\ran=
\n+\sum_{\bc\in\wt\cC}\frac{\o(\bc)}{|\bc|!}\,
t^{|\bc|}\cos\lan\t(\bc),k\ran,
\end{multline*}
where the series converge absolutely for all $t\in\C$. Thus, the
identity \er{hkfH}  has been proved. We rewrite this Fourier series
in the standard form
\begin{multline*}
\Tr e^{tH(k)}=\n+\sum_{\bc\in\wt\cC}\frac{\o(\bc)}{|\bc|!}\,
t^{|\bc|}\cos\lan\t(\bc),k\ran\\=\n+
\sum_{\mm\in\Z^d}\sum_{\bc\in\wt\cC^\mm}\frac{\o(\bc)}{|\bc|!}\,
t^{|\bc|}\cos\lan\mm,k\ran=\n+\sum\limits_{\mm\in\Z^d}\gh_\mm(t)\cos\lan\mm,k\ran,
\end{multline*}
where $\gh_\mm(t)$ are given in \er{THFS}. Integrating \er{hkfH}
over $k\in\T^d$  and using \er{IVFu} and \er{inex}, we obtain
\er{hkfHI}. \qq \BBox

\medskip

\no\textbf{Proof of Corollary \ref{CHC1}.} We have
\[\lb{reid}
\big(H(k)-\l I\big)^{-1}=-\frac1\l\,\bigg(I-\frac{H(k)}\l\bigg)^{-1}=
-\frac1\l\sum_{n=0}^\iy\frac{H^n(k)}{\l^n}\,,
\]
where the series converges absolutely for $|\l|>\|H(k)\|$. Taking
the trace of  \er{reid} and using  \er{TrnH} and \er{deTnk}, we
obtain
\begin{multline*}
\Tr\big(H(k)-\l I\big)^{-1}=-\frac1\l\sum_{n=0}^\iy\frac{\Tr H^n(k)}{\l^n}=
-\frac\n\l-\sum_{n=1}^\iy\frac{\cT_n(k)}{\l^{n+1}}\\=
-\frac\n\l-\sum_{n=1}^\iy\frac1{\l^{n+1}}\sum_{\bc\in\wt\cC_n}\o(\bc)
e^{-i\lan\t(\bc),k\ran}=
-\frac\n\l-\sum_{\bc\in\wt\cC}\frac{\o(\bc)}{\l^{|\bc|+1}}\,e^{-i\lan\t(\bc),k\ran}=
-\frac\n\l-\sum_{\bc\in\wt\cC}\frac{\o(\bc)}{\l^{|\bc|+1}}\,\cos\lan\t(\bc),k\ran\\=
-\frac\n\l-\sum_{\mm\in\Z^d}\sum_{\bc\in\wt\cC^\mm}
\frac{\o(\bc)}{\l^{|\bc|+1}}\,\cos\lan\mm,k\ran
=-\frac\n\l-\sum_{\mm\in\Z^d}\cR_\mm(\l)\cos\lan\mm,k\ran,
\end{multline*}
where $R_\mm(\l)$ are given in \er{reexHF}. Thus, \er{reexH0} and
\er{reexHF}  have been proved.

Integrating \er{reexH0} over $k\in\T^d$ and using \er{IVFu} and \er{inex}, we obtain \er{TrRI}. Since $\|H(k)\|\leq\|H\|$ for all $k\in\T^d$, the series in \er{TrRI} converges absolutely for $|\l|>\|H\|$. \qq \BBox

\medskip

Repeating a cycle $\bc$ $m$ times, we obtain \emph{$m$-multiple
$\bc^m$}  of $\bc$. If $\bc$ is not a $m$-multiple of a cycle with
$m\geq2$, $\bc$ is called \emph{prime}. Two cycles are
\emph{equivalent} if one is obtained from another by a cyclic
permutation of edges. All cyclic permutations of edges of a prime
cycle $\bc$ give $|\bc|$ distinct prime cycles having the same
length $|\bc|$, weight $\o(\bc)$ and index $\t(\bc)$. Note that for a non-prime cycle
$\bc$ this is not true anymore, since there are distinct cyclic
permutations of edges of $\bc$ giving the same cycles.

An equivalence class of a prime cycle $\bc$ will be denoted by
$\bc_*$. We define the length, index and weight of $\bc_*$ as the
length, index and weight of any representative of $\bc_*$:
\[\lb{ecpc}
|\bc_*|=|\bc|, \qqq \t(\bc_*)=\t(\bc), \qqq \o(\bc_*)=\o(\bc) \qqq \textrm{ for some \; $\bc\in\bc_*$ }.
\]
This definition is correct, since all equivalent cycles have the
same length, index and weight. We define \\
$\bu$  the set $\wt\cP$ of all equivalence classes of prime cycles in the modified fundamental graph $\wt\cG_*$;\\
$\bu$  the set $\wt\cP^0$ of all equivalence classes from $\wt\cP$ with zero index:
\[\lb{wtP0}
\wt\cP^0=\{\bc_*\in\wt\cP:\t(\bc_*)=0\}.
\]

We formulate the trace formulas for the resolvent and the determinant formulas for the fiber Schr\"odinger operators in terms of the equivalence classes of prime cycles in the modified fundamental graph.

\begin{corollary}
\lb{CHC2}
Let $H(k)$, $k\in\T^d$, be the fiber Schr\"odinger operator defined by \er{Hvt'} -- \er{fado}. Then

i) The trace of the resolvent of $H(k)$ has the following expansion
\[\lb{trfS}
\Tr\big(H(k)-\l I\big)^{-1}=-\frac\n\l
+\frac1\l\sum_{\bc_*\in\wt\cP}|\bc_*|\Big(1-\frac{\l^{|\bc_*|}}{\o(\bc_*)}\,
e^{i\lan\t(\bc_*),k\ran}\Big)^{-1}, \qqq \n=\#\cV_*,
\]
and
\[\lb{TrRI2}
\frac1{(2\pi)^d}\int_{\T^d}\Tr \big(H(k)-\l I\big)^{-1}dk
=-\frac\n\l+\frac1\l\sum_{\bc_*\in\wt\cP^0}|\bc_*|
\Big(1-\frac{\l^{|\bc_*|}}{\o(\bc_*)}\Big)^{-1}\,.
\]
The series in \er{trfS} converges absolutely for $|\l|>\|H(k)\|$, and the series in \er{TrRI2} converges absolutely for $|\l|>\|H\|$.

ii) The determinant of $I-tH(k)$ satisfies
\[\lb{TrloH}
\det\big(I-tH(k)\big)=\prod_{\bc_*\in\wt\cP}\big(1-
e^{-i\lan\t(\bc_*),k\ran}\o(\bc_*)t^{|\bc_*|}\big),
\]
and
\[\lb{TrloHI}
\frac1{(2\pi)^d}\int_{\T^d}\Tr\,\log\big(I-t H(k)\big)dk=\log\prod_{\bc_*\in\wt\cP^0}\big(1-
\o(\bc_*)t^{|\bc_*|}\big),
\]
where $|\bc_*|$, $\o(\bc_*)$, $\t(\bc_*)$ are defined in \er{ecpc}. The product in \er{TrloH} converges for $|t|<\|H(k)\|^{-1}$, and the product in \er{TrloHI} converges for $|t|<\|H\|^{-1}$.
\end{corollary}

\no\textbf{Proof.} \emph{i}) Due to \er{reexH0}, we have
\[\lb{prTr}
\Tr\big(H(k)-\l I\big)^{-1}=-\frac\n\l-
\frac1\l\sum_{\bc\in\wt\cC}\frac{\o(\bc)}{\l^{|\bc|}}\,e^{-i\lan\t(\bc),k\ran},
\]
where the series converges absolutely for $|\l|>\|H(k)\|$. Each cycle $\bc\in\wt\cC$ can be expressed in
a unique way in the form $\bc=m\bc_0$ for some $m\in\N$ and some
prime cycle $\bc_0$, and
$$
|\bc|=m|\bc_0|,\qqq \t(\bc)=m\t(\bc_0),\qqq \o(\bc)=\o^m(\bc_0).
$$
Then, using that each equivalence class $\bc_*$ of a prime cycle
$\bc$  consists of $|\bc|$ distinct prime cycles having the same
length $|\bc|$, index $\t(\bc)$ and weight $\o(\bc)$, we can rewrite
\er{prTr} in the form
\begin{multline*}
\Tr\big(H(k)-\l I\big)^{-1}=
-\frac\n\l-\frac1{\l}\sum_{\bc_*\in\wt\cP}|\bc_*|\sum_{m=1}^\iy
\bigg(\frac{\o(\bc_*)}{\l^{|\bc_*|}}\,e^{-i\lan\t(\bc_*),k\ran}\bigg)^m\\
=-\frac\n\l-\frac1\l\sum_{\bc_*\in\wt\cP}|\bc_*|
\bigg(\frac1{1-\frac{\o(\bc_*)}{\l^{|\bc_*|}}\,e^{-i\lan\t(\bc_*),k\ran}}-1\bigg)=
-\frac\n\l
+\frac1\l\sum_{\bc_*\in\wt\cP}|\bc_*|\Big(1-\frac{\l^{|\bc_*|}}{\o(\bc_*)}\,
e^{i\lan\t(\bc_*),k\ran}\Big)^{-1}.
\end{multline*}
Similarly, using \er{TrRI}, we obtain
\begin{multline*}
\frac1{(2\pi)^d}\int_{\T^d}\Tr\big(H(k)-\l I\big)^{-1}dk
=-\frac\n\l-\sum_{\bc\in\wt\cC^0}
\frac{\o(\bc)}{\l^{|\bc|+1}}=
-\frac\n\l-\frac1\l\sum_{\bc_*\in\wt\cP^0}|\bc_*|\sum_{m=1}^\iy
\bigg(\frac{\o(\bc_*)}{\l^{|\bc_*|}}\bigg)^m\\
=-\frac\n\l-\frac1\l\sum_{\bc_*\in\wt\cP^0}|\bc_*|
\bigg(\frac1{1-\frac{\o(\bc_*)}{\l^{|\bc_*|}}}-1\bigg)=
-\frac\n\l+\frac1\l\sum_{\bc_*\in\wt\cP^0}|\bc_*|
\Big(1-\frac{\l^{|\bc_*|}}{\o(\bc_*)}\Big)^{-1},
\end{multline*}
where the series converge absolutely for $|\l|>\|H\|$.

\emph{ii}) Let $|t|<\|H(k)\|^{-1}$. Then, using the power series for the matrix logarithm and the identities \er{TrnH}, \er{deTnk}, we have
\begin{multline}\lb{rLn}
\Tr\,\log\big(I-t H(k)\big)=
-\sum_{n=1}^\iy\frac{t^n}{n}\,\Tr H^n(k)
=-\sum_{n=1}^\iy\frac{t^n}{n}\,\sum_{\bc\in\wt\cC_n}e^{-i\lan\t(\bc),k\ran}\o(\bc)\\=
-\sum_{\bc\in\wt\cC}\frac{t^{|\bc|}}{|\bc|}\,e^{-i\lan\t(\bc),k\ran}\o(\bc)=
-\sum_{\bc_*\in\wt\cP}\sum_{m=1}^\iy\frac1{m}\,\Big(e^{-i\lan \t(\bc_*),
k\ran}\o(\bc_*)t^{|\bc_*|}\Big)^m
\\=\sum_{\bc_*\in\wt\cP}\log\Big(1-e^{-i\lan\t(\bc_*),k\ran}\o(\bc_*)t^{|\bc_*|}\Big)
=\log\prod_{\bc_*\in\wt\cP}\Big(1-e^{-i\lan\t(\bc_*),k\ran}\o(\bc_*)t^{|\bc_*|}\Big).
\end{multline}
This and the formula  $\det e^A=e^{\Tr A}$, where $A$ is any complex square
matrix, yield \er{TrloH}.

Integrating the first identity in \er{rLn} over $k\in\T^d$ and using \er{IVFu}, we obtain
\begin{multline*}
\frac1{(2\pi)^d}\int_{\T^d}\Tr\,\log\big(I-t H(k)\big)dk
=-\sum_{n=1}^\iy\frac{t^n}{n}\,\cT_{n,0}=
-\sum_{n=1}^\iy\frac{t^n}{n}\sum_{\bc\in\wt\cC_n^0}\o(\bc)
=-\sum_{\bc\in\wt\cC^0}\frac{t^{|\bc|}}{|\bc|}\,\o(\bc)\\=
-\sum_{\bc_*\in\wt\cP^0}\sum_{m=1}^\iy\frac1{m}\,\big(\o(\bc_*)t^{|\bc_*|}\big)^m
=\sum_{\bc_*\in\wt\cP^0}\log\big(1-\o(\bc_*)t^{|\bc_*|}\big)
=\log\prod_{\bc_*\in\wt\cP^0}\big(1-\o(\bc_*)t^{|\bc_*|}\big),
\end{multline*}
where the product converges for $|t|<\|H\|^{-1}$. \qq \BBox

\begin{remarks}
1) One can obtain the trace of $f\big(H(k)\big)$, where $f(x)$,
$|x|\leq\|H(k)\|$,  is any function defined via a convergent power
series.

2) The trace formulas for the heat kernel and for the resolvent and the
determinant formula for the operator $-\D(k)$, where $\D(k)=\vk-A(k)$ is the
fiber Laplacian, are given by the same identities \er{hkfH} --
\er{hkfHI}, \er{reexH0} -- \er{TrRI} and \er{trfS} -- \er{TrloHI},
where $\o(\bc)=\o(\bc,0)$ and $\o(\bc_*)=\o(\bc_*,0)$.

3) For a positive self-adjoint operator $T$ the $\G$-determinant is defined by
\[
\det_\G(T)=\exp\Tr_\G(\log T).
\]
where $\Tr_\G T$ is the $\G$-trace of $T$, see \er{GTr}. Thus, using \er{DIGT}, the identity \er{TrloHI} can be rewritten in the form
\[
\det_\G(I-tH)=\prod_{\bc_*\in\wt\cP^0}\big(1-
\o(\bc_*)t^{|\bc_*|}\big).
\]
This is an analog of the determinant formula for the adjacency operator on a regular periodic graph.

4) The \emph{Ihara zeta function} of a $\G$-periodic graph $\cG$ is defined by
\[
Z(u)=\prod_{\bc_*\in\cP^0}\big(1-u^{|\bc_*|}\big)^{-1}
\]
for all $u\in\C$ sufficiently small so that the infinite product converges.
Here $\cP^0$ denotes the set of equivalence classes of prime non-backtracking cycles with zero indices in the fundamental graph $\cG_*$.

5) The \emph{$L$-function}, a generalization of zeta functions, for the fiber adjacency operator $A(k)$ on the fundamental graph $\cG_*$ is defined by
\[
\cL(u,k)=
\prod_{\bc_*\in\cP}\Big(1-e^{-i\lan\t(\bc_*),k\ran}u^{|\bc_*|}\Big)^{-1},
\]
see \cite{AS87}. Here $\cP$ is the set of equivalence classes of prime cycles in the fundamental graph $\cG_*$. Thus, for the fiber adjacency $A(k)$ we can rewrite \er{TrloH} in the form
\[\lb{TrloD}
\det\big(I-u A(k)\big)=\cL(u,k)^{-1}.
\]
\end{remarks}

\section{Examples of trace formulas}
\setcounter{equation}{0}
\lb{Sec4}

\subsection{Examples} We prove Examples \ref{ExSL}, where we present the traces of the fiber Schr\"odinger and adjacency operators on the square lattice. We also consider the  Schr\"odinger operator with periodic potentials on the Kagome lattice.

\medskip

\no \textbf{Proof of Example \ref{ExSL}.} The fundamental graph $\dL_*^2=(\cV_*,\cA_*)$  of the square lattice $\dL^2$ with the periods $a_1$, $a_2$ (see Fig.~\ref{slex}\emph{a}) consists of 4 vertices $x_1,x_2,x_3,x_4$ and 8 edges
$$
\begin{array}{llll}
\be_1=(x_1,x_4), \qqq & \be_2=(x_4,x_2),\qqq & \be_3=(x_2,x_3),\qqq & \be_4=(x_3,x_1),\\
\be_5=(x_4,x_1), \qqq & \be_6=(x_2,x_4),\qqq & \be_7=(x_3,x_2),\qqq & \be_8=(x_1,x_3)
\end{array}
$$
with indices
$$
\begin{array}{llll}
\t(\be_1)=(0,0), \qqq & \t(\be_2)=(0,0), \qqq & \t(\be_3)=(0,0),\qqq & \t(\be_4)=(0,0)\\[2pt]
\t(\be_5)=(0,1), & \t(\be_6)=(1,0), & \t(\be_7)=(0,-1), & \t(\be_8)=(-1,0),
\end{array}
$$
and their inverse edges, see Fig.~\ref{slex}\emph{b}. Note that the
square lattice  $\dL^2$ and its fundamental graph $\dL^2_*$ are
bipartite. The fiber Schr\"odinger operator $H(k)$, defined by
\er{Hvt'} -- \er{fado} on $\dL_*^2$  in the standard orthonormal
basis of $\ell^2(\cV_*)$ has the form
\[\lb{efHSL}
H(k)=A(k)+v, \qqq\textrm{where}\qqq v=\diag\big(V_x-4\big)_{x\in\cV_*},
\]
and the $4\ts4$ matrix $A(k)$ of the fiber adjacency operator is given by
\[\lb{efHSL1}
A(k)=\left(
\begin{array}{cc}
\O_2 & B(k) \\
B^*(k) & \O_2\\
\end{array}\right),\qq
B(k)=\left(
\begin{array}{ll}
1+e^{ik_1} & 1+e^{ik_2}\\
1+e^{-ik_2} & 1+e^{-ik_1}  \\
\end{array}\right),\qq k=(k_1,k_2)\in\T^2,
\]
where $\O_2$ is the $2\ts 2$ zero matrix.

We calculate the traces of $A^n(k)$ and $H^n(k)$, $n=1,2,3$, using
formulas  \er{TrAo}, \er{TrH123} and \er{todT}.  There are no cycles
of length one (i.e., loops) in the fundamental graph $\dL_*^2$.
Then, by \er{TrAo} as $n=1$ and the first identity in \er{TrH123},
we have
\[\lb{Tr2_13}
\begin{aligned}
&\Tr A(k)=\sum_{\bc\in\cC_1}\cos\lan\t(\bc),k\ran=0,\\
 &H(k)=\sum_{x\in\cV_*}\big(V_x-4\big)=-16+\sum_{s=1}^4V_s.
\end{aligned}
\]
Of course, these identities also follow from the explicit form
\er{efHSL}, \er{efHSL1} of the fiber operators $H(k)$ and $A(k)$.

Since there are no multiple edges in the square lattice $\dL^2$,
then,  by \er{cN0g1}, the number of all cycles of length 2 with zero
index $\cN_2^0=\#\cA_*=16$. The graph $\dL^2_*$ also has the
following proper cycles (i.e., cycles without backtracking) of
length 2:
$$
\bc_{2,1}=(\be_1,\be_5), \qqq \bc_{2,2}=(\be_2,\be_6), \qqq \bc_{2,3}
=(\be_3,\be_7),\qqq \bc_{2,4}=(\be_4,\be_8),
$$
with indices
$$
\t(\bc_{2,1})=(0,1),\qqq \t(\bc_{2,2})=(1,0),\qqq \t(\bc_{2,3})=(0,-1),
\qqq \t(\bc_{2,4})=(-1,0),
$$
their cyclic permutations (two permutations for each cycle
$\bc_{2,s}$,  $s=1,2,3,4$) and their reverse cycles. Thus, by
\er{TrAo}  as $n=2$ and the second identity in \er{TrH123}, we have
$$
\Tr A^2(k)=\sum_{\bc\in\cC_2}\cos\lan\t(\bc),k\ran=8\cF(k),\qqq \cF(k)=2+\cos k_1+\cos k_2,
$$
$$
\Tr H^2(k)=\sum_{s=1}^4v_s^2+\sum_{\bc\in\cC_2}\cos\lan\t(\bc),k\ran=
\sum_{s=1}^4v_s^2+\Tr A^2(k), \qqq v_s=V_s-4.
$$

Since $\dL^2_*$ is bipartite, then there are no cycles of length 3
in $\dL_*^2$.  Then, by \er{TrAo} as $n=3$ and \er{todT}, we obtain
$$
\Tr A^3(k)=\sum_{\bc\in\cC_3}\cos\lan\t(\bc),k\ran=0,
$$
$$
\Tr H^3(k)=\sum_{s=1}^4v_s^3+\frac32\sum_{\bc\in\cC_2}v(\bc)\cos\lan\t(\bc),k\ran
=\sum_{s=1}^4\big(v_s^3+6\cF(k)v_s\big). \qqq \BBox
$$

\begin{remark}
The square lattice $\dL^2$ (Fig.\ref{slex}\emph{a}) is bipartite,
but its minimal  fundamental graph $\dL^2/\G_0$, where $\G_0$ is a
lattice generated by the orthonormal basis $e_1,e_2$ of $\R^2$,
consists of one vertex and two loop edges at this vertex, and is
non-bipartite. But if we consider the sublattice $\G\ss\G_0$ with
periods $a_1=2e_1$ and $a_2=2e_2$, the fundamental graph
$\dL^2_*=\dL^2/\G$ of $\dL^2$ is bipartite (see
Fig.\ref{slex}\emph{b}).
\end{remark}

\setlength{\unitlength}{1.0mm}
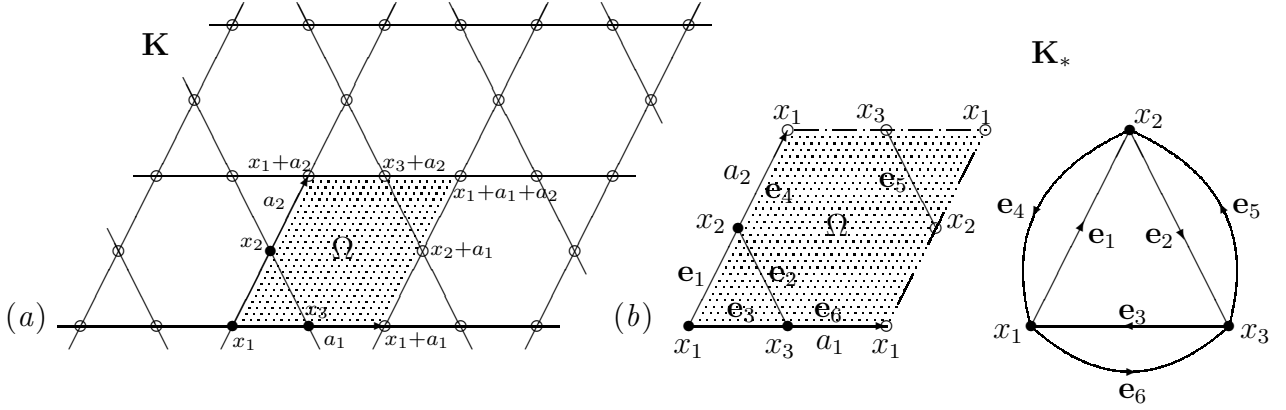
\begin{figure}[h]
\centering
\unitlength 1mm 
\linethickness{0.4pt}
\ifx\plotpoint\undefined\newsavebox{\plotpoint}\fi 

\begin{picture}(160,55)(0,0)
\bezier{25}(21,15)(26,25)(31,35)
\bezier{25}(22,15)(27,25)(32,35)
\bezier{25}(23,15)(28,25)(33,35)
\bezier{25}(24,15)(29,25)(34,35)
\bezier{25}(25,15)(30,25)(35,35)
\bezier{25}(26,15)(31,25)(36,35)
\bezier{25}(27,15)(32,25)(37,35)
\bezier{25}(28,15)(33,25)(38,35)
\bezier{25}(29,15)(34,25)(39,35)
\bezier{25}(30,15)(35,25)(40,35)
\bezier{25}(31,15)(36,25)(41,35)
\bezier{25}(32,15)(37,25)(42,35)
\bezier{25}(33,15)(38,25)(43,35)
\bezier{25}(34,15)(39,25)(44,35)
\bezier{25}(35,15)(40,25)(45,35)
\bezier{25}(36,15)(41,25)(46,35)
\bezier{25}(37,15)(42,25)(47,35)
\bezier{25}(38,15)(43,25)(48,35)
\bezier{25}(39,15)(44,25)(49,35)

\put(8.0,51){$\bK$}
\put(-3,15){\line(1,0){66.0}}
\put(7,35){\line(1,0){66.0}}
\put(17,55){\line(1,0){66.0}}
\put(-1.5,12){\line(1,2){23.0}}
\put(18.5,12){\line(1,2){23.0}}
\put(38.5,12){\line(1,2){23.0}}
\put(58.5,12){\line(1,2){23.0}}

\put(11.5,12){\line(-1,2){8.0}}
\put(31.5,12){\line(-1,2){18.0}}
\put(51.5,12){\line(-1,2){23.0}}
\put(66.5,22){\line(-1,2){18.0}}
\put(76.5,42){\line(-1,2){8.0}}

\put(20,15){\vector(1,0){20.0}}
\put(20,15){\vector(1,2){10.0}}

\put(0,15){\circle{1.5}}
\put(10,15){\circle{1.5}}
\put(20,15){\circle*{1.5}}
\put(30,15){\circle*{1.5}}
\put(40,15){\circle{1.5}}
\put(50,15){\circle{1.5}}
\put(60,15){\circle{1.5}}

\put(5,25){\circle{1.5}}
\put(25,25){\circle*{1.5}}
\put(45,25){\circle{1.5}}
\put(65,25){\circle{1.5}}

\put(10,35){\circle{1.5}}
\put(20,35){\circle{1.5}}
\put(30,35){\circle{1.5}}
\put(40,35){\circle{1.5}}
\put(50,35){\circle{1.5}}
\put(60,35){\circle{1.5}}
\put(70,35){\circle{1.5}}

\put(15,45){\circle{1.5}}
\put(35,45){\circle{1.5}}
\put(55,45){\circle{1.5}}
\put(75,45){\circle{1.5}}

\put(20,55){\circle{1.5}}
\put(30,55){\circle{1.5}}
\put(40,55){\circle{1.5}}
\put(50,55){\circle{1.5}}
\put(60,55){\circle{1.5}}
\put(70,55){\circle{1.5}}
\put(80,55){\circle{1.5}}

\put(32.0,12.5){$\scriptstyle a_1$}
\put(24.0,31.0){$\scriptstyle a_2$}
\put(20,12.2){$\scriptstyle x_1$}
\put(22.0,36.0){$\scriptstyle x_1+a_2$}
\put(39.7,36.0){$\scriptstyle x_3+a_2$}
\put(40.0,12.5){$\scriptstyle x_1+a_1$}
\put(21,25.5){$\scriptstyle x_2$}
\put(46.0,24.5){$\scriptstyle x_2+a_1$}
\put(49.0,32.0){$\scriptstyle x_1+a_1+a_2$}
\put(29.5,16.5){$\scriptstyle x_3$}
\put(33,24){$\Omega$}
\put(-10,15.0){(\emph{a})}

\put(125.0,50){$\bK_*$}
\put(98,27){$\Omega$}
\multiput(94,41)(5,0){5}{\line(1,0){3}}
\bezier{30}(107.5,18)(108.5,20)(109.5,22)
\bezier{30}(111.5,26)(112.5,28)(113.5,30)
\bezier{30}(115.5,34)(116.5,36)(117.5,38)

\put(96.5,16.3){$\be_6$}
\put(85.0,16.3){$\be_3$}
\put(90.5,21){$\be_2$}
\put(78.5,21){$\be_1$}
\put(90.0,32){$\be_4$}
\put(105,33.5){$\be_5$}

\put(80,15){\vector(1,0){26.0}}
\put(80,15){\vector(1,2){13.0}}

\bezier{30}(106,15)(112.5,28)(119,41)

\bezier{30}(81,15)(87.5,28)(94,41)
\bezier{30}(82,15)(88.5,28)(95,41)
\bezier{30}(83,15)(89.5,28)(96,41)
\bezier{30}(84,15)(90.5,28)(97,41)
\bezier{30}(85,15)(91.5,28)(98,41)
\bezier{30}(86,15)(92.5,28)(99,41)
\bezier{30}(87,15)(93.5,28)(100,41)
\bezier{30}(88,15)(94.5,28)(101,41)
\bezier{30}(89,15)(95.5,28)(102,41)
\bezier{30}(90,15)(96.5,28)(103,41)
\bezier{30}(91,15)(97.5,28)(104,41)
\bezier{30}(92,15)(98.5,28)(105,41)
\bezier{30}(93,15)(99.5,28)(106,41)
\bezier{30}(94,15)(100.5,28)(107,41)
\bezier{30}(95,15)(101.5,28)(108,41)
\bezier{30}(96,15)(102.5,28)(109,41)
\bezier{30}(97,15)(103.5,28)(110,41)
\bezier{30}(98,15)(104.5,28)(111,41)
\bezier{30}(99,15)(105.5,28)(112,41)
\bezier{30}(100,15)(106.5,28)(113,41)
\bezier{30}(101,15)(107.5,28)(114,41)
\bezier{30}(102,15)(108.5,28)(115,41)
\bezier{30}(103,15)(109.5,28)(116,41)
\bezier{30}(104,15)(110.5,28)(117,41)
\bezier{30}(105,15)(111.5,28)(118,41)

\put(80,15){\circle*{1.5}}
\put(93,15){\circle*{1.5}}
\put(106,15){\circle{1.5}}

\put(93,15){\line(-1,2){6.5}}
\put(112.5,28){\line(-1,2){6.5}}

\put(86.5,28){\circle*{1.5}}
\put(112.5,28){\circle{1.5}}

\put(93,41){\circle{1.5}}
\put(106,41){\circle{1.5}}
\put(119,41){\circle{1.5}}

\put(96.7,12.0){$a_1$}
\put(84.5,34.5){$a_2$}
\put(78.0,11.5){$x_1$}
\put(91.0,42.5){$x_1$}
\put(102.0,42.7){$x_3$}
\put(104.0,11.5){$x_1$}
\put(80.9,28.0){$x_2$}
\put(113.9,28.0){$x_2$}
\put(116.0,42.5){$x_1$}
\put(90.0,11.5){$x_3$}
\put(70,15.0){(\emph{b})}

\put(125,15){\line(1,0){26.0}}
\put(125,15){\line(1,2){13.0}}
\put(151,15){\line(-1,2){13.0}}
\put(125,15){\circle*{1.5}}
\put(120.0,13.5){$x_1$}
\put(151,15){\circle*{1.5}}
\put(152.5,13.5){$x_3$}
\put(138,41){\circle*{1.5}}
\put(138.5,42){$x_2$}
\bezier{200}(125,15)(138,3)(151,15)
\bezier{200}(125,15)(120,33)(138,41)
\bezier{200}(151,15)(156,33)(138,41)

\put(136.5,16){$\be_3$}
\put(136.5,5.5){$\be_6$}
\put(132.5,26){$\be_1$}
\put(140,26){$\be_2$}
\put(120.5,30){$\be_4$}
\put(151.5,30){$\be_5$}

\put(131.0,27.0){\vector(1,2){1.0}}
\put(144.2,28.5){\vector(1,-2){1.0}}
\put(138,15.0){\vector(-1,0){1.0}}
\put(150.8,29.5){\vector(-1,2){1.0}}
\put(126.2,31.5){\vector(-1,-2){1.0}}
\put(138,9.0){\vector(1,0){1.0}}
\end{picture}
\vspace{-0.8cm} \caption{ \emph{a}) Kagome lattice $\bK$; the
fundamental cell $\Omega$  is shaded; \quad \emph{b})
the~fundamental graph $\bK_*$.} \lb{FEx2}
\end{figure}

\begin{example}\lb{EKL}
Let $H=A-\vk+V$ be the Schr\"odinger operator on the Kagome
lattice $\bK$.  We assume that the potential $V$ has the periods
$a_1$, $a_2$, see Fig.~\ref{FEx2}a. Then the traces of $A^n(k)$ and
$H^n(k)$, $k=(k_1,k_2)\in\T^2$, given by \er{TrAo} and \er{TrnH},
for $n=1,2,3$ have the form
\[\lb{TAKL}
\Tr A(k)=0,\qqq
\Tr A^2(k)=12+4\mF(k),\qqq
\Tr A^3(k)=12+12\mF(k),
\]
\[\lb{Tr2_1'}
\begin{array}{l}
\Tr H(k)=-12+\sum\limits_{s=1}^3V_s,\qqq
\Tr H^2(k)=\Tr A^2(k)+\sum\limits_{s=1}^3v_s^2,\\[4pt]
\Tr H^3(k)=\Tr
A^3(k)+\sum\limits_{s=1}^3v_s^3+12\sum\limits_{s=1}^3v_s
\\[2pt]
\hspace{17mm}+6\,(v_1+v_3)\cos k_1+6\,(v_1+v_2)\cos
k_2+6\,(v_2+v_3)\cos(k_1-k_2),
\end{array}
\]
where
$$
\mF(k)=\cos k_1+\cos k_2+\cos(k_1-k_2), \qqq v_s=V_s-4, \qqq V_s=V_{x_s},\qqq s=1,2,3.
$$
\end{example}

\no \textbf{Proof.} The fundamental graph $\bK_*$ of the Kagome
lattice $\bK$  consists of three vertices $x_1,x_2,x_3$, six edges
$$
\be_1=(x_1,x_2),\qq \be_2=(x_2,x_3),\qq \be_3=(x_3,x_1),\qq
\be_4=(x_2,x_1),\qq \be_5=(x_3,x_2),\qq \be_6=(x_1,x_3)
$$
with indices
$$
\begin{array}{lll}
\t(\be_1)=(0,0), \qqq & \t(\be_2)=(0,0), \qqq & \t(\be_3)=(0,0),\\[2pt]
\t(\be_4)=(0,1), & \t(\be_5)=(1,-1), & \t(\be_6)=(-1,0),
\end{array}
$$
and their inverse edges, see Fig.~\ref{FEx2}\emph{b}. Then the fiber
Schr\"odinger  operator $H(k)$ defined by \er{Hvt'} -- \er{fado} on
$\bK_*=(\cV_*,\cA_*)$  in the standard orthonormal basis of
$\ell^2(\cV_*)$ has the form
\[\lb{efHk}
H(k)=A(k)+v, \qqq\textrm{where}\qqq v=\diag\big(V_x-4\big)_{x\in\cV_*},
\]
and the $3\ts3$ matrix $A(k)$ of the fiber adjacency operator is given by
\[\lb{effL}
A(k)=\left(
\begin{array}{ccc}
0 & 1+e^{ik_2} & 1+e^{ik_1} \\
1+e^{-ik_2} & 0 & 1+e^{i(k_1-k_2)}\\
1+e^{-ik_1} & 1+e^{-i(k_1-k_2)}  & 0
\end{array}\right),\qqq k=(k_1,k_2)\in\T^2.
\]

We calculate the traces of $A^n(k)$ and $H^n(k)$, $n=1,2,3$, using
formulas  \er{TrAo}, \er{TrH123} and \er{todT}. There are no cycles
of length one (i.e., loops) in the fundamental graph $\bK_*$. Then,
by \er{TrAo} as $n=1$ and the first identity in \er{TrH123}, we have
\[\lb{Tr2_1}
\begin{aligned}
&\Tr A(k)=\sum_{\bc\in\cC_1}\cos\lan\t(\bc),k\ran=0,
\\
&H(k)=\sum_{x\in\cV_*}\big(V_x-4\big)=-12+\sum_{s=1}^3V_s.
\end{aligned}
\]

Since there are no multiple edges in the Kagome lattice $\bK$,
then, by \er{cN0g1}, the number of all cycles of length 2 with zero
index $\cN_2^0=\#\cA_*=12$. The graph $\bK_*$ also has the following
proper cycles (i.e., cycles without backtracking) of length 2:
$$
\bc_{2,1}=(\be_6,\be_3), \qqq \bc_{2,2}=(\be_1,\be_4), \qqq \bc_{2,3}=(\be_2,\be_5),
$$
with indices
$$
\t(\bc_{2,1})=(-1,0),\qqq \t(\bc_{2,2})=(0,1),\qqq \t(\bc_{2,3})=(1,-1),
$$
their cyclic permutations (two permutations for each cycle
$\bc_{2,s}$,  $s=1,2,3$) and their reverse cycles. Thus, by
\er{TrAo} as $n=2$ and the second identity in \er{TrH123}, we have
$$
\Tr A^2(k)=\sum\limits_{\bc\in\cC_2}\cos\lan\t(\bc),k\ran=12+
4\cos k_1+4\cos k_2+4\cos(k_1-k_2)=12+4\mF(k),
$$
$$
\Tr H^2(k)=\sum_{s=1}^3v_s^2+\sum\limits_{\bc\in\cC_2}\cos\lan\t(\bc),k\ran=
\sum_{s=1}^3v_s^2+\Tr A^2(k), \qqq v_s=V_s-4.
$$
There are no cycles of length 3 with backtracking in $\bK_*$, since
there are no loops and each backtracking contributes 2 in the cycle
length. Thus, all cycles of length 3 are proper cycles
$$
\begin{array}{llll}
\bc_{3,1}=(\be_1,\be_2,\be_3),\qq & \bc_{3,2}=(\be_1,\be_2,\ul\be_6),\qq &
\bc_{3,3}=(\be_2,\be_3,\ul\be_4),\qq & \bc_{3,4}=(\be_3,\be_1,\ul\be_5),\\
\bc_{3,5}=(\be_6,\be_5,\be_4),  & \bc_{3,6}=(\be_6,\be_5,\ul\be_1),
& \bc_{3,7}=(\be_5,\be_4,\ul\be_3), & \bc_{3,8}=(\be_4,\be_6,\ul\be_2)
\end{array}
$$
with indices
$$
\begin{array}{llll}
\t(\bc_{3,1})=(0,0), \qq & \t(\bc_{3,2})=(1,0), \qq &
\t(\bc_{3,3})=(0,-1), \qq & \t(\bc_{3,4})=(-1,1),\\ \t(\bc_{3,5})=(0,0),
\qq & \t(\bc_{3,6})=(0,-1), \qq & \t(\bc_{3,7})=(1,0), \qq & \t(\bc_{3,8})=(-1,1),
\end{array}
$$
their cyclic permutations (three permutations for each cycle
$\bc_{3,s}$,  $s\in\N_8$) and their reverse cycles. Then, by
\er{TrAo} as $n=3$ and \er{todT}, we obtain
$$
\Tr A^3(k)=\sum\limits_{\bc\in\cC_3}\cos\lan\t(\bc),k\ran=
12+12\cos k_1+12\cos k_2+12\cos(k_1-k_2)=12+12\mF(k),
$$
$$
\begin{aligned}
&\Tr H^3(k)=\sum_{s=1}^3v_s^3
+\frac32\sum\limits_{\bc\in\cC_2}v(\bc)\cos\lan\t(\bc),k\ran
+\sum\limits_{\bc\in\cC_3}\cos\lan\t(\bc),k\ran=\sum_{s=1}^3v_s^3
+12\sum\limits_{s=1}^3v_s\\
&+6\,(v_1+v_3)\cos k_1+6\,(v_1+v_2)\cos k_2+6\,(v_2+v_3)\cos(k_1-k_2)+\Tr A^3(k).\qqq \BBox
\end{aligned}
$$

\section{Trace formulas for normalized Laplacians}
\setcounter{equation}{0}
\lb{Sec5}

\subsection{Normalized Laplacians on periodic graphs}
In this section we discuss trace formulas for normalized Laplacians
on periodic graphs $\cG=(\cV,\cA)$. We define the \emph{normalized} Laplacian
$\D_\gn$ on $\ell^2(\cV)$ by
\begin{equation}\lb{DNLA}
(\D_\gn f)_x=\sum_{(x,y)\in\cA}\frac{f_y}{\sqrt{\vk_x\vk_y}}\,,
 \qqq f\in\ell^2(\cV),\qqq x\in\cV,
\end{equation}
where $\vk_x$ is the degree of the vertex $x$. The sum in \er{DNLA} is
taken  over all edges from $\cA$ starting at the vertex $x$. In the
literature the operator $\D_\gn$ is usually called \emph{the transition
operator} and the normalized Laplacian is defined by $I-\D_\gn$, where
$I$ is the identity operator. But since $-\D_\gn$ and $I-\D_\gn$ differ only
by a shift we will call $\D_\gn$ the normalized Laplacian.

It is known (see, e.g., \cite{MW89}) that the normalized Laplacian
$\D_\gn$  is  a bounded self-adjoint operator on $\ell^2(\cV)$ and
its spectrum $\s(\D_\gn)$ is a closed subset of the segment
$[-1,1]$, containing the point 1, i.e.,
\begin{equation}
\lb{mp}
1\in\s(\D_\gn)\subseteq[-1,1].
\end{equation}

\begin{remarks} 1) If $\cG$ is a regular graph of degree $\vk_+$, then
the normalized Laplacian $\D_\gn$, the combinatorial Laplacian
$\D$ and the adjacency operator $A$ defined in \er{Sh}, \er{ALO}
(and, consequently, their spectra) are related by the simple
identities
$$
A=\vk_+\D_\gn,\qqq \D=\vk_+I-A=\vk_+(I-\D_\gn).
$$
However, in the case of an arbitrary
graph the spectra of $\D$ and $\D_\gn$, in spite of many similar
properties, may have significant differences.

2) Recall that when we consider the combinatorial Laplacians,
without loss of generality we may assume that there are no loops in
the periodic graph. But for the normalized Laplacian this is not
true anymore.
\end{remarks}

We recall Proposition 1.1 from \cite{KS18}.

\begin{theorem}\label{TFD1}
Let the Hilbert space $\mH$ be defined by \er{Hisp}. Then the
normalized  Laplacian $\D_\gn$ on $\ell^2(\cV)$ given by \er{DNLA}
has the following decomposition into a constant fiber direct
integral
\[
\lb{razN}
\begin{aligned}
& U\D_\gn U^{-1}=\int^\oplus_{\T^d}\D_\gn(k){dk\/(2\pi)^d}\,,
\end{aligned}
\]
where $U:\ell^2(\cV)\to\mH$ is some unitary operator (the Gelfand
transform), and the fiber Laplacian $\D_\gn(k)$, $k\in\T^d$, on
$\ell^2(\cV_*)$ is given by
\[
\label{ftro} \big(\D_\gn(k)f\big)_x=\sum_{\be=(x,\,y)\in\cA_*}
\frac{e^{i\lan\t(\be), \,k\ran}}{\sqrt{\vk_x\vk_y}}\,f_y, \qqq
f\in\ell^2(\cV_*),\qqq x\in \cV_*,
\]
where $\vk_x$ is the degree of the vertex $x$, and $\t(\be)$ is  the
index of the edge $\be\in\cA_*$ defined by \er{in}, \er{dco}.
\end{theorem}

Each fiber operator $\D_\gn(k)$, $k\in\T^{d}$, has $\n$ real
eigenvalues $\m_j(k)$, $j\in\N_\n$, $\n=\#\cV_*$,  which are labeled
in non-decreasing order (counting multiplicities) by
\[
\label{eq.3} \m_{1}(k)\leq\m_{2}(k)\leq\ldots\leq\m_{\nu}(k),
\qqq \forall\,k\in\T^{d}.
\]
Since $\D_\gn(k)$ is self-adjoint and analytic in $k\in\T^{d}$, each
$\m_j(\cdot)$, $j\in\N_\n$, is a real and piecewise analytic
function on the torus $\T^{d}$ and creates the \emph{spectral band}
$\s_j(\D_\gn)$ given by
\[
\lb{ban.1}
\s_j(\D_\gn)=[\m_j^-,\m_j^+]=\m_j(\T^{d}).
\]
Note that $\m_\n^+=\m_\n(0)=1$. Then the spectrum of the Laplacian
$\D_\gn$  on the periodic graph $\cG$ is given by
\[\lb{spec}
\s(\D_\gn)=\bigcup_{k\in\T^d}\s\big(\D_\gn(k)\big)=
\bigcup_{j=1}^{\nu}\s_j(\D_\gn)=\s_{ac}(\D_\gn)\cup \s_{fb}(\D_\gn),
\]
where $\s_{ac}(\D_\gn)$ is the absolutely continuous spectrum, which is a
union of non-degenerate bands, and $\s_{fb}(\D_\gn)$ is the set of
all flat bands (eigenvalues of infinite multiplicity).

\subsection{Trace formulas for normalized Laplacians} Recall that the cycle sets $\cC$, $\cC_n$, $\cC^\mm$, $\cC_n^\mm$ and the numbers $\cN_n$ and $\cN_n^\mm$ are defined in Subsection \ref{SoC}.

For each cycle $\bc=(\be_1,\ldots,\be_n)\in\cC$ we define the \emph{weight}
\[\lb{cywe}
\o_\gn(\bc)=\frac1{\vk_{x_1}\ldots\vk_{x_n}}\,, \qq \textrm{where}\qq
 \be_s=(x_{s},x_{s+1})\in\cA_*, \qq s\in\N_n, \qq x_{n+1}=x_1,
\]
and $\vk_x$ is the degree of the vertex $x$.

\medskip

Similarly to the functionals $\cT_n(k)$ given by \er{deTnk} we
define  the functionals $T_n(k)$ as the finite sum over cycles form
$\cC_n$ given by
\[
\lb{deTnkN}
T_n(k)=\sum_{\bc\in\cC_n}\o_\gn(\bc)e^{-i\lan\t(\bc),k\ran} =
\sum_{\bc\in\cC_n}\o_\gn(\bc)\cos\lan\t(\bc),k\ran,
\]
where we have used that $\o_\gn(\bc)=\o_\gn(\ul\bc\,)$ and
$\t(\bc)=-\t(\ul\bc\,)$  for each $\bc\in\cC_n$.

We formulate some properties of the functionals $T_n(k)$. Let
$\vk_-$ and $\vk_+$  be the minimum and maximum vertex degrees, i.e.
\[\lb{vl+-}
\vk_-=\min\limits_{v\in\cV_*}\vk_v,\qqq
\vk_+=\max\limits_{v\in\cV_*}\vk_v.
\]

\begin{proposition}\lb{PpofN}
Let the functional $T_n(k)$ be defined by \er{deTnkN}. Then

i) $T_n(k)$  satisfies
\[\lb{esgh}
\big|T_n(k)\big|\leq T_n(0)\leq
\frac{\cN_n}{\vk_-^n}\,, \qqq \big|T_n(k)\big|\leq\n,
\]
where $\vk_-$ is given in \er{vl+-}.

ii) $T_n(k)$ has the following finite Fourier series
\[\lb{FsTrHN}
T_n(k)=\sum_{\mm\in\Z^d\atop\|\mm\|\leq n\t_+}
T_{n,\mm}\cos\lan\mm,k\ran,\qqq
0\leq\frac{\cN_n^m}{\vk_+^n}\leq T_{n,\mm}=\sum_{\bc\in\cC_n^\mm}\o_\gn(\bc)\leq
\frac{\cN_n^m}{\vk_-^n}\,,
\]
and
\[\lb{IVFuN}
\frac1{(2\pi)^d}\int_{\T^d}T_n(k) dk=T_{n,0},\qqq
0\leq\frac{\cN_n^0}{\vk_+^n}\leq T_{n,0}=\sum_{\bc\in\cC_n^0}\o_\gn(\bc)\leq
\frac{\cN_n^0}{\vk_-^n}\,,
\]
where $\t_+$ and $\o_\gn(\bc)$ are defined in \er{cNe0} and
\er{cywe}, respectively,  and $\|\cdot\|$ is the standard norm in
$\R^d$.

iii) If $\cG_*$ is a regular graph of degree $\vk_+$, then
\[\label{ghkRR}
T_n(k)=\frac1{\vk_+^n}\sum_{\bc\in\cC_n}\cos\lan\t(\bc),k\ran=
\frac1{\vk_+^n}\sum_{\mm\in\Z^d\atop\|\mm\|\leq n\t_+}\cN_n^\mm\cos\lan\mm,k\ran,
\]
\[\lb{IVFuNR}
\frac1{(2\pi)^d}\int_{\T^d}T_n(k) dk=\frac{\cN_n^0}{\vk_+^n}\,.
\]
\end{proposition}

\no \textbf{Proof.} The proof of items \emph{i}), \emph{ii}) is
similar to  the proof of Proposition \ref{Ppof}.

\emph{iii}) Let $\cG_*$ be regular of degree $\vk_+$.  Then
$\o_\gn(\bc)=\frac1{\vk^{|\bc|}_+}$ for each cycle $\bc\in\cC$, and
the identities \er{deTnkN}, \er{FsTrHN}, \er{IVFuN}  have the form
\er{ghkRR}, \er{IVFuNR}. \qq $\BBox$

\medskip

Now we formulate trace formulas for the fiber normalized Laplacians.

\begin{theorem}\lb{TFNL0} Let $\D_\gn(k)$, $k\in\T^d$, be the
fiber normalized Laplacian defined by \er{ftro} on the fundamental
graph $\cG_*=(\cV_*,\cA_*)$. Then the following statements hold
true.

i) For each $n\in\N$ the trace of $\D_\gn^n(k)$ satisfies
\[\lb{TrnN0}
\Tr \D_\gn^n(k)=\sum_{j=1}^\n\m_j^n(k)=T_n(k), \qqq
\frac1{(2\pi)^d}\int_{\T^d}\Tr\D_\gn^n(k)dk=T_{n,0},\qqq \n=\#\cV_*,
\]
where $T_n(k)$ and $T_{n,0}$ are defined in \er{deTnkN} and \er{IVFuN}, respectively.

ii) The fundamental graph $\cG_*$ is bipartite iff
\[\lb{KbfgN}
\Tr \D_\gn^n(0)=0\qqq \textrm{for all odd $n\leq\n$}.
\]

iii) The periodic graph $\cG$ is bipartite iff
\[\lb{KbpgN}
\int_{\T^d}\Tr \D_\gn^n(k)dk=0\qqq \textrm{for all odd $n$}.
\]
Moreover, the condition \er{KbpgN} can not be reduced, i.e., for any $s\in\N$ there exists a non-bipartite periodic graph $\cG$ such that $\int_{\T^d}\Tr \D_\gn^n(k)dk=0$ for all odd $n<s$.
\end{theorem}

We omit the proof since it repeats the proof of Theorems \ref{TPG},
\ref{TPG2} and Theorem \ref{TFao}.\emph{iii})-\emph{iv}).

\subsection{Trace formulas for the heat kernel, the resolvent and the determinant formula}

We present trace formulas for the heat kernel and the resolvent of the fiber normalized Laplacian.

\begin{corollary} \lb{CHFO}
Let $\D_\gn(k)$, $k\in\T^d$, be the fiber Laplacian defined by
\er{ftro} on the fundamental graph $\cG_*=(\cV_*,\cA_*)$. Then the
trace of  the heat kernel $e^{t\D_\gn(k)}$, $t\in\C$, satisfies:
\[\lb{THFL}
\Tr e^{t\D_\gn(k)}=\n+\sum_{n=1}^\iy\frac{t^n}{n!}\,T_n(k)=
\n+\sum_{\bc\in\cC}\frac{\o_\gn(\bc)}{|\bc|!}\,
t^{|\bc|}\cos\lan\t(\bc),k\ran,\qqq \n=\#\cV_*,
\]
\[\lb{TDFS}
\Tr
e^{t\D_\gn(k)}=\n+\sum_{\mm\in\Z^d}h_\mm(t)\cos\lan\mm,k\ran,\qqq
h_\mm(t)=\sum_{\bc\in\cC^\mm}\frac{\o_\gn(\bc)}{|\bc|!}\,
t^{|\bc|},
\]
\[\lb{hkfLI}
\frac1{(2\pi)^d}\int_{\T^d}\Tr e^{t\D_\gn(k)}dk=\n+\sum_{n=1}^\iy\frac{t^n}{n!}\,T_{n,0}=
\n+\sum_{\bc\in\cC^0}\frac{\o_\gn(\bc)}{|\bc|!}\, t^{|\bc|},
\]
where $T_n(k)$ and $T_{n,0}$ are given in \er{deTnkN} and
\er{IVFuN},  respectively. The series in \er{THFL} -- \er{hkfLI}
converge absolutely for all $t\in\C$.
\end{corollary}

\begin{corollary} \lb{Tr1} Let $\D_\gn(k)$, $k\in\T^d$, be
 the fiber Laplacian defined by \er{ftro} on the fundamental graph $\cG_*=(\cV_*,\cA_*)$.
 Then the trace of the resolvent of $\D_\gn(k)$ has the following expansions:
\[\lb{reexN0}
\Tr\big(\D_\gn(k)-\l I\big)^{-1}=
-\frac\n\l-\sum_{n=1}^\iy\frac{T_n(k)}{\l^{n+1}}=
-\frac\n\l-\sum_{\bc\in\cC}\frac{\o_\gn(\bc)}{\l^{|\bc|+1}}\,\cos\lan\t(\bc),k\ran,
\]
\[\lb{rexN0}
\Tr\big(\D_\gn(k)-\l I\big)^{-1}=-\frac\n\l-\sum_{\mm\in\Z^d}R_\mm(\l)\cos\lan\mm,k\ran,\qq
\textrm{where}\qq R_\mm(\l)=\sum_{\bc\in\cC^\mm}
\frac{\o_\gn(\bc)}{\l^{|\bc|+1}}\,,
\]
\[\lb{TrRIN}
\frac1{(2\pi)^d}\int_{\T^d}\Tr \big(\D_\gn(k)-\l I\big)^{-1}dk
=-\frac\n\l-\sum_{n=1}^\iy\frac{T_{n,0}}{\l^{n+1}}
=-\frac\n\l-\sum_{\bc\in\cC^0}\frac{\o_\gn(\bc)}{\l^{|\bc|+1}},
\]
where $T_n(k)$ and $T_{n,0}$ are given in \er{deTnkN} and
\er{IVFuN}, respectively. The series
in \er{reexN0} -- \er{TrRIN} converge absolutely for $|\l|>1$.
\end{corollary}

The proof of Corollaries \ref{CHFO} and \ref{Tr1} is
omitted as it is similar to the proof of Corollaries \ref{CHC},
and \ref{CHC1}.

\medskip

We introduce \\
$\bu$  the set $\cP$ of all equivalence classes of prime cycles in the  fundamental graph $\cG_*$;\\
$\bu$  the set $\cP^0$ of all equivalence classes from $\cP$ with zero index:
\[\lb{NLP0}
\cP^0=\{\bc_*\in\cP:\t(\bc_*)=0\}.
\]
For an equivalence class  $\bc_*\in\cP$ of a prime cycle $\bc\in\cC$ we define the weight $\o_\gn(\bc_*)$ as the weight of any representative of
$\bc_*$:
\[\lb{ecpcWN}
\o_\gn(\bc_*)=\o_\gn(\bc) \qq \textrm{ for some \; $\bc\in\bc_*$}.
\]

We formulate the trace formulas for the resolvent and the determinant formulas for the fiber normalized Laplacian in terms of the equivalence classes of prime cycles in the fundamental graph.

\begin{corollary} \lb{CHFO1}
Let $\D_\gn(k)$, $k\in\T^d$, be the fiber Laplacian defined by \er{ftro}. Then

i) The trace of the resolvent of $\D_\gn(k)$ has the following expansion
\[\lb{reexN}
\Tr\big(\D_\gn(k)-\l I\big)^{-1}=-\frac\n\l
+\frac1\l\sum_{\bc_*\in\cP}|\bc_*|\Big(1-\frac{\l^{|\bc_*|}}{\o_\gn(\bc_*)}\,
e^{i\lan\t(\bc_*),k\ran}\Big)^{-1}, \qqq \n=\#\cV_*,
\]
and
\[\lb{TrRI2N}
\frac1{(2\pi)^d}\int_{\T^d}\Tr \big(\D_\gn(k)-\l I\big)^{-1}dk
=-\frac\n\l+\frac1\l\sum_{\bc_*\in\cP^0}|\bc_*|
\Big(1-\frac{\l^{|\bc_*|}}{\o_\gn(\bc_*)}\Big)^{-1}\,.
\]
The series in \er{reexN}, \er{TrRI2N} converge absolutely for $|\l|>1$.

ii) The determinant of $I-t\D_\gn(k)$ satisfies
\[\lb{Trlo}
\det\big(I-t\D_\gn(k)\big)=
\prod_{\bc_*\in\cP}\big(1-e^{-i\lan\t(\bc_*),k\ran}\o_\gn(\bc_*)t^{|\bc_*|}\big),
\]
and
\[\lb{TrloDI}
\frac1{(2\pi)^d}\int_{\T^d}\Tr\log\big(I-t\D_\gn(k)\big)dk=
\log\prod_{\bc_*\in\cP^0}\big(1-
\o_\gn(\bc_*)t^{|\bc_*|}\big),
\]
where $|\bc_*|$, $\t(\bc_*)$ and $\o_\gn(\bc_*)$ are defined in \er{ecpc} and \er{ecpcWN}. The products in \er{Trlo}, \er{TrloDI} converge for $|t|<1$.
\end{corollary}

The proof of Corollary \ref{CHFO1} is omitted as it is similar to the proof of Corollary \ref{CHC2}.

\medskip

\no\textbf{Acknowledgments.} \footnotesize Our study was supported
by the   RFBR grant No. 19-01-00094.


\medskip

\begin{thebibliography}{9999}
\setlength{\itemsep}{-\parskip}\footnotesize
\bibitem[AS87]{AS87} Adachi, T.; Sunada, T. Twisted Perron-Frobenius theorem and $L$-functions, J. Funct. Anal. 71 (1987), 1--46.

\bibitem[Ah87]{Ah87} Ahumada, G. Fonctions periodiques et formule
 des traces de Selberg sur les arbres, C. R. Acad. Sci. Paris 305(1987), 709--712.

\bibitem[A76]{A76} Atiyah, M.F. Elliptic operators, discrete groups and von Neumann algebras, Asterisque 32--33 (1976), 43--72.

\bibitem[Br91]{Br91} Brooks, R. The Spectral Geometry of $k$-Regular Graphs,
J. d'Analyse, 57 (1991), 120--151.

\bibitem[CJK15]{CJK15} Chinta, G.; Jorgenson, J.; Karlsson, A. Heat kernels on regular graphs and generalized Ihara zeta function formulas, Monatsh. Math. 178 (2015), no. 2, 171--190.

\bibitem[GIL08]{GIL08} Guido, D.; Isola, T.; Lapidus M. L. Ihara's zeta function for periodic graphs and its approximation in the amenable case, Journal of Functional Analysis 255 (2008), 1339--1361.

\bibitem[GIL09]{GIL09} Guido, D.; Isola, T.; Lapidus, M.L. A trace on fractal graphs and the Ihara Zeta function, Trans. Amer. Math. Soc. 361 (2009), 3041--3070.

\bibitem[HS99]{HS99} Higuchi, Y.; Shirai, T. The spectrum
of magnetic Schr\"odinger operators on a graph with periodic
structure, J. Funct. Anal., 169 (1999), 456--480.

\bibitem [IK12]{IK12} Isozaki, H.; Korotyaev, E. Inverse problems, trace
formulae for discrete Schr\"odinger operators, Ann. Henri
Poincar\'e, 13 (2012), 751--788.

\bibitem [K17]{K17}  Korotyaev, E. Trace formulae for Schr\"odinger operators on lattice, preprint arXiv:1702.01388.

\bibitem [K21]{K21} Korotyaev, E. Trace formulas for time periodic complex Hamiltonians on lattice, preprint arXiv:2101.03370.

\bibitem [KL18]{KL18} Korotyaev, E.;  Laptev, A. Trace formulae for
Schr\"odinger operators with complex-valued potentials on cubic
lattices,  Bulletin of Mathematical Sciences,  8 (2018),
453--475.

\bibitem[KS14]{KS14} Korotyaev, E.; Saburova, N.
Schr\"odinger operators on periodic discrete graphs, J. Math. Anal. Appl.,
 420 (2014), no. 1, 576--611.

 \bibitem[KS17]{KS17}  Korotyaev, E.; Saburova, N.
Magnetic Schr\"odinger operators on periodic discrete graphs, J.
Funct. Anal., 272 (2017), 1625--1660.

\bibitem[KS18]{KS18} Korotyaev, E.; Saburova, N. Spectral estimates
 for Schr\"odinger operators on periodic discrete graphs, St. Petersburg Math. J. 30 (2018),
  no.~4, 667--698.

\bibitem[KS20]{KS20} Korotyaev, E.; Saburova, N. Invariants for Laplacians
 on periodic graphs, Math. Ann. 377 (2020), 723--758.

\bibitem[KS21a]{KS21a} Korotyaev, E.; Saburova, N. Two-sided estimates of total bandwidth for Schr\"odinger operators on periodic graphs, preprint arXiv:2106.08661.

\bibitem[KS21b]{KS21b} Korotyaev, E.; Saburova, N. On continuous spectrum
 of magnetic Schr\"odinger operators on periodic discrete graphs, preprint arXiv:2101.05571.

\bibitem[KoS03]{KoS03} Kotani, M.; Sunada, T. Spectral geometry of crystal lattices, Contemporary Math.,  338 (2003), 271--305.

\bibitem[LPS19]{LPS19} Lenz, D.; Pogorzelski, F; Schmidt, M. The Ihara zeta function for infinite graphs, Trans. Amer. Math. Soc. 371 (2019), 5687--5729.

\bibitem [MN15]{MN15} Malamud, M.; Neidhardt, H. Trace formulas for additive
 and non-additive perturbations. Adv. Math. 274 (2015), 736--832.

\bibitem[Mn07]{Mn07} Mn\"{e}v, P. Discrete Path Integral Approach
 to the Selberg Trace Formula for Regular Graphs, Commun. Math. Phys. 274 (2007), 233--241.

\bibitem[MW89]{MW89} Mohar, B.; Woess, W. A survey on spectra
 of infinite graphs, Bull. London Math. Soc. 21 (1989), no.~3, 209--234.

\bibitem[NG04]{NG04} Novoselov, K.S.; Geim, A.K. et al,
Electric field effect in atomically thin carbon films, Science 22 October,
306 (2004), no. 5696, 666--669.

\bibitem[OGS09]{OGS09} Oren, I.; Godel, A.; Smilansky, U. Trace formulae
 and spectral statistics for discrete Laplacians on regular graphs. I. J. Phys. A, 42 (2009), 415101.

\bibitem[S86]{S86} Sunada, T. $L$-functions in geometry and some applications, Springer Lecture Notes in Math. 1201, 1986, pp.266--284.

\bibitem[SS92]{SS92} Sy, P.W.; Sunada, T. Discrete Schr\"odinger
 operator on a graph, Nagoya Math. J., 125 (1992), 141--150.

\end{thebibliography}
\end{document}